\documentclass[preprint]{imsart}
\RequirePackage[OT1]{fontenc}
\RequirePackage{amsthm,amsmath}
\RequirePackage[numbers]{natbib}
\RequirePackage[colorlinks,citecolor=blue,urlcolor=blue]{hyperref}

\usepackage[utf8]{inputenc}
\usepackage{amsfonts, amsmath, amssymb, amsthm, bbm, color, enumerate, graphicx, mathtools, tikz, hyperref, relsize}
\usepackage[top=1in, bottom=1in, left=1in, right=1in]{geometry}

\definecolor{expcol}{rgb}{1.0,0.5,0.5}
\definecolor{ecol}{rgb}{0.0, 0.0, 1.0}

\newcommand{\ZZ}{\mathbf{Z}}
	
\newcommand{\yy}{\mathbf{y}}			

\newcommand{\xx}{\mathbf{x}}
\newcommand{\kk}{\mathbf{k}}
\newcommand{\XX}{\mathbf{X}}
\newcommand{\YY}{\mathbf{Y}}

\newcommand{\one}{\mathbf{1}}

\newcommand{\clu}{{\mathcal{U}}}
\newcommand{\clv}{{\mathcal{V}}}

\newcommand{\NN}{\mathbb{N}}
\newcommand{\Om}{\Omega}
\newcommand{\clf}{{\mathcal{F}}}
\newcommand{\cls}{{\mathcal{S}}}
\newcommand{\clp}{{\mathcal{P}}}
\newcommand{\Pmb}{{\mathbb{P}}}

\newcommand{\floor}[1]{\lfloor #1 \rfloor}

\newcommand{\R}{\mathbb{R}}
\newcommand{\Z}{\mathbf{Z}}

\newcommand{\U}{\mathbf{U}}
\newcommand{\V}{\mathbf{V}}

\newcommand{\RR}{\mathbb{R}}
\newcommand{\cld}{\mathcal{D}}
\newcommand{\clc}{\mathcal{C}}
\newcommand{\cll}{\mathcal{L}}

\newcommand{\blow}{\underline{b}}
\newcommand{\rinv}{R^{-1}}
\newcommand{\sigmalow}{\underline{\sigma}}
\newcommand{\sigmahigh}{\overline{\sigma}}
\newcommand{\rinvk}[1]{\left(R|_{#1}\right)^{-1}}
\newcommand{\muk}[1]{\mu|_{#1}}
\newcommand{\bk}[1]{\mathbf{b}^{(#1)}}
\newcommand{\ak}[1]{a^{(#1)}}

\newcommand{\blowk}[1]{\blow^{(#1)}}
\newcommand{\bnorm}[1]{\|#1 \|_{1, \beta}}
\newcommand{\supnorm}[1]{\|#1\|_{\infty}}

\newcommand{\gbd}{{\boldsymbol{g}}}

\newcommand{\pbd}{{\boldsymbol{p}}}

\newcommand{\xbd}{{\boldsymbol{x}}}
\newcommand{\ybd}{{\boldsymbol{y}}}\newcommand{\Ybd}{{\boldsymbol{Y}}}
\newcommand{\zbd}{{\boldsymbol{z}}}\newcommand{\Zbd}{{\boldsymbol{Z}}}

\newcommand{\Ybf}{{\mathbf{Y}}}

\newtheorem{theorem}{Theorem}
\newtheorem{corollary}[theorem]{Corollary}
\newtheorem{lemma}[theorem]{Lemma}

\newtheorem{remark}{Remark}
\newtheorem{definition}{Definition}

\newtheorem{property}{Property}

\let\plainqed\qedsymbol
\newcommand{\claimqed}{$\lrcorner$}

\hypersetup{
 colorlinks=true,
 linkcolor=blue,          
 citecolor=red,       
 filecolor=red,   
 urlcolor=red, 
 pdftitle={},
 pdfauthor={},
 pdfsubject={},
 pdfkeywords={},
 linktocpage=true
}

\usepackage{fancyhdr}
 
\usepackage{natbib}

\usepackage{amssymb}
\usepackage{amsmath}
\usepackage{amsthm}
\newcommand{\Expect}[1]{\operatorname{\mathbb{E}}\left[#1\right]}

\newcommand{{\LPC}}{\textbf{LPC}}

\newcommand{\Reals}{\mathbb{R}}
\newcommand{\x}{\mathbf{x}}
\newcommand{\y}{\mathbf{y}}
\newcommand{\z}{\mathbf{z}}
\newcommand{\D}{\mathbf{D}}
\newcommand{\0}{\mathbf{0}}
\newcommand{\bb}{\mathbf{b}}
\newcommand{\RBM}{\operatorname{RBM}}

\newtheorem{thm}{Theorem}

\newtheorem{rem}[thm]{Remark}

\begin{document}

\begin{frontmatter}

\title{Long Time Behavior of Finite and Infinite Dimensional Reflected Brownian Motions
}
\runtitle{Ergodicity of RBM}
\author{\fnms{Sayan} \snm{Banerjee}\ead[label=e1]{sayan@email.unc.edu}}
\and
\author{\fnms{Amarjit} \snm{Budhiraja}\ead[label=e2]{amarjit@unc.edu}}
\affiliation{University of North Carolina, Chapel Hill}
\runauthor{Banerjee and Budhiraja}
\address{Department of Statistics \\and Operations Research\\
University of North Carolina\\
Chapel Hill, NC 27599\\
In Celebration of Prof. Rajeeva Karandikar's 65th Birthday
}

\begin{abstract}
This article presents a review of some old and new results on the long time behavior of reflected diffusions. First, we present a summary of prior results on construction, ergodicity and geometric ergodicity of reflected diffusions in the positive orthant $\mathbb{R}^d_+$, $d \in \mathbb{N}$. The geometric ergodicity results, although very general, usually give implicit convergence rates due to abstract couplings and Lyapunov functions used in obtaining them. This leads us to some recent results on an important subclass of reflected Brownian motions (RBM) (constant drift and diffusion coefficients and oblique reflection at boundaries), known as the Harrison-Reiman class, where explicit rates of convergence are obtained as functions of the system parameters and underlying dimension. In addition,  sufficient conditions on system parameters of the RBM are provided under which local convergence to stationarity holds at a `dimension-free' rate, that is, for any fixed $k \in \mathbb{N}$, the rate of convergence of the $k$-marginal to equilibrium does not depend on the dimension of the whole system. Finally, we study the long time behavior of infinite dimensional rank-based diffusions, including the well-studied infinite Atlas model. The gaps between the ordered particles evolve as infinite dimensional RBM and this gap process has uncountably many explicit product form stationary distributions. Sufficient conditions for initial configurations to lie in the weak domain of attraction of the various stationary distributions are provided. Finally, it is shown that, under conditions, all of these explicit stationary distributions are extremal (equivalently, ergodic) and, in some sense, the only product form invariant probability distributions. Proof techniques involve a pathwise analysis of RBM using explicit synchronous and mirror couplings and constructing Lyapunov functions.
\end{abstract}


\end{frontmatter}

\section{Introduction}
In this article, we will provide a survey of some results on reflected Brownian motions (RBM) in polyhedral domains that focus on the long-time behavior of such processes. We will consider the settings of both finite and infinite dimensional reflected diffusions. 

The simplest such model is a one dimensional RBM $\{X(t)\}_{t\ge 0}$ defined as
\begin{equation}\label{eq:zeq}
	X(t) = x + B(t) - \inf_{0\le s \le t} \left\{(x+B(s))\wedge 0\right\}, t \ge 0, \; x \in \RR_+,
\end{equation}
where $B$ is a standard Brownian motion given on some probability space. Occasionally we will write $X(t) = X(t;x)$ to emphasize dependence on the initial condition. There is an equivalent way to describe this process using the formulation of the Skorohod problem
\cite{skorokhod1961stochastic}.

The Skorohod problem on the positive real line $\RR_+$ is to find for a given $\yy \in \cld([0,\infty):\RR)$ (throughout, for a Polish space $S$, $\cld([0,\infty):S)$ will denote the space of functions from $[0,\infty)$ to $S$ that are right continuous and have left limits(RCLL), equipped with the usual Skorohod topology), with $\yy(0)\ge 0$, a pair of paths $\xx,\kk \in \cld([0,\infty):\RR)$
such that for all $t\ge 0$, $\xx(t) = \yy(t)+\kk(t) \ge 0$, $\kk(0)=0$, $\kk$ is nondecresing, and $\kk$ increases only at instants $t$ when $\xx(t)=0$, namely
$$\int_{[0,\infty)} 1_{\{\xx(s)>0\}} d\kk(s) =0.$$
The above problem has a unique solution for every $\yy \in \cld([0,\infty):\RR)$ and we write $\xx = \Gamma_1(\yy)$. Also, it is easily verified that $X$ defined by \eqref{eq:zeq} can equivalently be written as $X = \Gamma_1(x+ B)$. The formulation of a Skorohod problem also allows one to define a general reflected diffusion that on $(0,\infty)$ behaves as a diffusion with the infinitesimal generator (evaluated on smooth test functions $f$)
$$\cll f(x) = \frac{1}{2}\sigma^2(x) f''(x) + b(x) f'(x),\; x\in (0,\infty)$$
for suitable coefficient functions $\sigma$ and $b$. Such a reflected diffusion can be constructed by solving the stochastic equation
\begin{equation}\label{eq:genrefdif}
	X(t) = \Gamma_1\left(x + \int_0^{\cdot} b(X(s)) ds + \int_0^{\cdot} \sigma(X(s)) dB(s)\right)(t), \;  t \ge 0
\end{equation}
in a pathwise fashion (e.g. when $b,\sigma$ are Lipschitz functions) (see \cite{andore}). In the special case when $\sigma(z) \equiv \sigma >0$ and $b(z) = b \in \RR$, $X$ is referred to as a RBM with drift $b$ and diffusion coefficient $\sigma$, or simply a $(b,\sigma)$-RBM.
It is well known that a $(b,\sigma)$-RBM is positive recurrent if and only if $b<0$ and the unique stationary distribution is an exponential distribution with rate $\lambda = 2|b|/\sigma^2$. In the general setting of \eqref{eq:genrefdif}, under conditions on the coefficients (e.g., $b$ and $\sigma$ Lipschitz, $\sigma$ bounded, and $\sigma \sigma^T$ uniformly nondegenerate), it can be shown that if $\sup_{z\ge 0} b(z) <0$ then the diffusion is positive recurrent, although this time an explicit form stationary distribution is not available. One can also show that the law at time $t$ converges to the stationary distribution at a geometric rate in the total variation distance.

In this work we are interested in the higher dimensional analogues of the above process. 
We will only consider reflected Brownian motions, or more generally reflected jump-diffusions, in a finite or infinite nonnegative orthant, namely the state space will be $\RR_+^d$ where $d \in \NN \cup \{\infty\}$.
Such multi-dimensional RBM arise from problems in stochastic networks \cite{reiman_jackson_net,harrison1987brownian}, applications in mathematical finance \cite{fernholz2009stochastic}, scaling limits of interacting particle systems \cite{karatzas_skewatlas}, and in rank-based diffusion models for competing particles  \cite{karatzas_skewatlas,sarantsev}.  The next three sections will focus on the finite dimensional case while the last two sections will consider infinite dimensional settings.

\section{Finite dimensional RBM}
Consider for now $d<\infty$. Roughly speaking, a reflected Brownian motion in the nonnegative orthant $\RR_+^d$ behaves like an ordinary Brownian motion (with drift) in the interior of the orthant and when it reaches the boundary of the domain it is instantaneously pushed back, with the minimal force needed to keep it inside the domain, applied in some pre-specified directions of constraint. As in the one dimensional case, a natural way to give a rigorous formulation of such an object is through a Skorohod problem.

Fix vectors $r^1, \ldots, r^d \in \RR^d$. Here the vector $r^i$ gives the direction in which the state of the RBM is pushed when it is about to exit the domain $\RR_+^d$ from its $i$-th face $F_i = \{x \in \RR_+^d: x_i =0\}$. We refer to the $d\times d$ matrix $R = [r^1, \ldots , r^d]$ as the reflection matrix.
For each $x \in \partial \R_+^d = \cup_{i=1}^d F_i$, the set $r(x)$ of directions of constraints is defined as
\begin{equation}
	r(x) := \left\{ Rq: q = (q_1, \ldots , q_d)' \in \RR_+^d, \; q\cdot \one =1, \mbox{ and } q_i>0 \mbox{ only if } x_i=0\right\},
\end{equation}
where $\one = (1, 1, \ldots)'$.
The set $r(x)$ represents the collection of directions available for pushing the state of the process at the boundary point $x$. Note that if $x$ has exactly one coordinate $0$, the set $r(x)$ is a singleton.

The Skorohod problem associated with the domain $\RR_+^d$ and the reflection matrix $R$ is defined as follows. 
Given $\psi \in \cld([0,\infty): \RR^d)$ such that $\psi(0)\in \RR_+^d$, we say a pair of trajectories
 $\phi, \eta \in \cld([0,\infty): \RR^d)$ solve the {\bf Skorohod problem} for $\psi$ with respect to
  reflection matrix $R$, if and only if $\eta(0)=0$ and for all $t\ge 0$: (i) $\phi(t) = \psi(t) + \eta(t)$; (ii) $\phi(t) \in \RR_+^d$; 
  (iii) The total variation (with respect to the Euclidean norm on $\RR^d$) of $\eta$ on $[0,t]$, denoted as $|\eta|(t)$ is finite; (iv) $|\eta|(t) = \int_{[0,t]} 1_{\{\phi(s) \in \partial \RR_+^d\}} d|\eta|(s)$; 
   (v) there is a Borel measurable function $\gamma: [0,\infty) \to \RR^d$ such that $\gamma(t) \in r(\phi(t))$, $d|\eta|$-a.e., and $\eta(t) = \int_{[0,t]} \gamma(s) d|\eta|(s)$.
   The trajectory $\phi$ can be viewed as the reflected version of $\psi$ that stays in $\RR_+^d$ at all times. Property (iv) says that the pushing term $\eta$ is activated only when the path is about to exit the domain while property (v) ensures that the push is applied in the permissible directions of constraint, as specified by the reflection matrix $R$.
   
   On the domain $D \subset D_0 := \{\phi \in \cld([0,\infty): \RR^d): \phi(0) \in \RR_+^d\}$ on which there is a unique solution to the Skorohod problem we define the {\bf Skorohod map} (SM) $\Gamma$ as $\Gamma(\psi) = \phi$, if $(\phi, \psi-\phi)$ is the unique solution of the Skorohod problem.
   
   In many problems of interest it turns out that the Skorohod map is quite well behaved in the  sense that the following property is satisfied.
  \begin{property}\label{prop:regsp}
  	The Skorohod map is well defined on all of $D_0$ and the SM is Lipschitz continuous in the following sense: There exists a $K \in (0, \infty)$ such that for all $\phi_1, \phi_2 \in D_0$
	$$\sup_{0\le t <\infty} \|\Gamma(\phi_1)(t)- \Gamma(\phi_2)(t)\| \le K\sup_{0\le t <\infty} \|\phi_1(t)- \phi_2(t)\| .$$
  \end{property}
We refer the reader to \cite{dupish, dupuram} for sufficient conditions under which the above property holds. 
One important class of examples that arise naturally in many stochastic network models \cite{reiman_jackson_net,harrison1987brownian} and also  in the study of gaps between $d+1$ competing particles in rank-based diffusions (e.g. \cite{karatzas_skewatlas,sarantsev}), and where the above property is satisfied, is the so called {\bf Harrison-Reiman} class \cite{harrison1981reflected}. 
This corresponds to the family of models for which the  matrix $P := I - R^T$ is substochastic (non-negative entries and row sums bounded above by 1) and transient ($P^n \rightarrow 0$ as $n \rightarrow \infty$). This family will be the focus of much of our discussion in subsequent sections.

When Property \ref{prop:regsp} holds,  one can solve pathwise, as in the one dimensional case,  stochastic equations of the form
\begin{equation}\label{eq:genrefdif2}
	X(t) = \Gamma\left(x + \int_0^{\cdot} b(X(s)) ds + \int_0^{\cdot} \sigma(X(s)) dB(s)\right)(t), \;  t \ge 0,
\end{equation}
where $b: \RR_+^d \to \RR^d$ and $\sigma: \RR_+^d \to \RR^{d\times k}$ are suitable coefficient functions (e.g. Lipschitz continuous) and
$B$ is a standard $k$ dimensional Brownian motion on some probability space (cf. \cite{dupish}).
The solution in the case when $b(x) \equiv \mu$ and $\sigma(x) \equiv \D$, where $\mu \in \RR^d$ and $\D$ is a $d\times d$ matrix, will be of particular interest here and  will be referred to as the $(\mu, \D)$-RBM with reflection matrix $R$. In later sections, we will also use the abbreviation $\RBM(\mu,\Sigma,R)$, where $\Sigma = \D\D^T$.

The solution of \eqref{eq:genrefdif2} can be written in a somewhat more explicit form using boundary local times. For example, when $X(t;x)$ is a $(\mu, \D)$-RBM with reflection matrix $R$ starting from $x$, it can be represented as  the solution of
the equation
\begin{equation}\label{RBMdef}
X(t;\x) = \x + \D B(t) + \mu t + RL(t),
\end{equation}
 where $L$, referred to as  the local time process, is a $d$-dimensional non-decreasing continuous process satisfying 
 \begin{equation}\label{loctim}
 L(0)=0, \;\; \int_0^t X_i(s;\x)dL_i(s) = 0 \mbox{ for all } t>0 \mbox{ and } 1 \le i \le d.
 \end{equation}
 
 It turns out that there are important applications arising from stochastic network theory (see e.g. \cite{williams1998diffusion}) that lead to well-defined Brownian motions constrained in the positive orthant according to a certain reflection mechanism, but the associated Skorohod problem is not well-posed for a typical path in $\cld([0,\infty): \RR_+^d)$ (in fact uniqueness of solution fails even for linear paths). In order to cover such situations the papers \cite{taywil, reiwil} introduce and study {\bf semimartingale reflecting Brownian motions} (SRBM).
 
Suppose that $\Sigma$ is a $d\times d$ positive definite matrix. For $\x \in \RR_+^d$, an SRBM associated with the 
data $\left(\mu, \D, R\right)$ that starts from $\x$ is a continuous, 
$\left\{\mathcal{F}_{t}\right\}$-adapted $d$-dimensional process $X$, defined on some filtered probability 
space $\left(\Omega, \mathcal{F},\left\{\mathcal{F}_{t}\right\}_{t \geq 0}, P\right)$ such that: (i) 
$X(t)=\x + \D B(t)+\mu t+R L(t) \in \RR_+^d$ for all $t \geq 0, P$-a.s.,
(ii) $B$ is a $d$-dimensional standard $\left\{\mathcal{F}_{t}\right\}$-Brownian motion, (iii) $L$ is an $\left\{\mathcal{F}_{t}\right\}$-adapted $d$-dimensional process such that $L_{i}(0)=0$ for $i=1, \ldots, d, P$-a.s. For each $i=1, \ldots, d, L_{i}$ is continuous, nondecreasing and  $\int_{0}^{t} 1_{\left\{X_{i}(s) \neq 0\right\}} d L_{i}(s)=0$ for all $t \geq 0, P$-a.s.

A key condition in this context is that the matrix  $R$ is {\bf completely- $\mathcal{S}$} (cf. \cite{reiwil}), namely, for every $k \times k$ principal submatrix $G$ of $R$, there is a $k$-dimensional vector $v_{G}$ such that $v_{G} \geq 0$ and $R v_{G}>0$. Roughly speaking, the conditions says that at each point $x$ on the boundary there is a permissible direction of reflection, namely a $r \in r(x)$, that points inwards.
In \cite[Theorem 2]{reiwil}, it was shown that a necessary condition for the existence of an SRBM is that the reflection matrix $R$ is a completely- $\mathcal{S}$ matrix. The paper \cite{taywil} shows that the condition is sufficient as well for (weak) existence and uniqueness of SRBM to hold. Under this condition, the collection  $\{P_{\x}\}_{\x \in \RR_+^d}$,
where $P_{\x}$ is the law of the SRBM starting from $\x$, forms a $\RR_+^d$-valued strong Markov process. 

\subsection{Ergodicity.}
The first basic result on long-time behavior of RBM in $\RR_+^d$ is due to \cite{harrison1987brownian} which treats the Harrison-Reiman class of RBM. Suppose that the following holds ((A1) is simply a restatement of the Harrison-Reiman class.)

 \textbf{Assumptions: }
 \begin{itemize}
 \item[(A1)] The matrix $P := I - R^T$ is substochastic  and transient.  
 \item[(A2)] $\bb := -R^{-1}\mu > 0$.
 \item[(A3)] The matrix $\Sigma = \D\D^T$ is positive definite.
 \end{itemize}
It was shown in \cite{harrison1987brownian}  that, when the above assumptions (A1)-(A3) hold, then the 
$(\mu, \D)$-RBM is positive recurrent and consequently has a unique stationary distribution. Assumption (A2) is the well known `stability condition' which is sufficient for the existence of a stationary measure \cite[Section 6]{harrison1987brownian}. The condition is almost necessary in that if $\bb_i<0$ for some $i$ then the RBM is transient \cite{buddupsimp}. The matrix $\Sigma = DD^T$ gives the covariance matrix associated with the driving diffusion of \eqref{RBMdef} and (A3) guarantees that this diffusion process is elliptic. Ellipticity gives the uniqueness of the stationary distribution \cite{harrison1987brownian}. Under additional conditions on $R$ and $\D$, the paper \cite{harrisonwilliams_expo} shows that  the stationary distribution takes an explicit and  simple form given as a product of exponential distributions, however in general there is little one can say about this stationary distribution other than some properties of its tail (see e.g. the next section and \cite{budlee}).

The paper \cite{ABD} studied long time behavior of constrained reflected diffusions of the form in \eqref{eq:genrefdif2}. The basic   assumptions are that Property \ref{prop:regsp} holds and the coefficients $b, \sigma$ satisfy the following condition.\\
 
\noindent  \textbf{Assumptions: }
 \begin{itemize}
 \item[(D1)] $b$ and $\sigma$ are Lipschitz maps.
 \item[(D2)] $\sigma$ is bounded.
 \item[(D3)] The matrix $\sigma \sigma^T$ is uniformly nondegenerate.
 \end{itemize}
 In addition one needs a suitable stability condition.
The key stability condition in  \cite{ABD} is formulated in terms of the cone
$$
\mathcal{C} :=\left\{-\sum_{i=1}^{d} \alpha_{i} r^{i}: \alpha_{i} \geq 0, i \in\{1, \ldots, d\}\right\}.
$$
This cone $\mathcal{C}$ was introduced in \cite{buddupsimp} to characterize stability properties of deterministic paths arising from certain law of large number limits, under the long-time scaling, of a $(\mu, \D)$-RBM for which the associated reflection matrix $R$ satisfies Property \ref{prop:regsp}. The main result there showed that if $\mu \in \mathcal{C}^o$ then these deterministic paths are all attracted to the origin (namely the paths converge to $0$ as time becomes large) whereas if $\mu \in \clc^c$ then these paths diverge to $\infty$. Finally, when $\mu \in \partial \clc$,  all solutions of the deterministic model are bounded, and  for at least one initial condition the corresponding trajectory is not attracted to the origin. This stability characterization was instrumental in the study of the stability of constrained reflected diffusions in  \cite{ABD}. Let,
for $\delta \in(0, \infty)$, 
$$
\mathcal{C}(\delta) :=\{v \in \mathcal{C}: \operatorname{dist}(v, \partial \mathcal{C}) \geq \delta\} .
$$
The key stability assumption on the constrained diffusion model \eqref{eq:genrefdif2} made in \cite{ABD} is the following.\\

\noindent  \textbf{Assumption ($\mathcal{C}$-delta): }
  There exists a $\delta \in(0, \infty)$ and a bounded set $A \subseteq \RR_+^d$ such that for all $x \in \RR_+^d \backslash A$, $b(x) \in \mathcal{C}(\delta)$.

Under this stability assumption, together with assumptions (D1)-(D3) noted previously, \cite[Theorem 2.2]{ABD} shows that the Markov process  described by the constrained diffusion model \eqref{eq:genrefdif2} is positive recurrent and has a unique stationary distribution. The main proof idea there is to study stability properties of a certain deterministic collection of controlled reflected paths of the form
$$\psi_v(t;\x) \doteq \Gamma(\x+ \int_0^{\cdot} v(s) ds)(t), \; t \ge0$$
where $\{v(t)\}_{t\ge 0}$ is a locally integrable function that take values in $\clc(\delta)$ for all $t\ge 0$. 
Denote by $T(\x)$ the supremum of the hitting times to zero, taken over all such controlled paths, starting from $\x$. Then, 
the analysis of \cite{buddupsimp}  gives an explicit bound on $T(\cdot)$ as, $T(\x) \le C \|\x\|/\delta$ for some $C \in (0,\infty)$ and all $\x \in \RR_+^d$. The function $T(\cdot)$ turns out to be a natural Lyapunov function for the problem, although it lacks the required $C^2$-property for implementing the usual It\^{o}'s formula based methods for Lyapunov functions for diffusion processes. Nevertheless, the function is Lipschitz and that is enough to deduce the required uniform in time moment estimates needed to prove the main result.

The cone $\clc$ can also be used to give conditions for  stability  of constrained L\'{e}vy processes and more general constrained jump-diffusions in $\RR_+^d$ associated with a reflection matrix $R$ for which Property \ref{prop:regsp} holds. Such a jump-diffusion is given as a solution of an equation of the form
\begin{equation}\label{eq:genrefdif3}
	X(t) = \Gamma\left(x + \int_0^{\cdot} b(X(s)) ds + \int_0^{\cdot} \sigma(X(s)) dB(s)
	+ J(\cdot)\right)(t), \;  t \ge 0,
\end{equation}
where
$$J(t) = \int_{[0,t]\times \RR^n} h(\psi(X(s-), z)) [N(ds, dz)- \lambda(dz) ds] +
\int_{[0,t]\times \RR^n} h'(\psi(X(s-), z)) N(ds, dz),$$
$N$ is a Poisson random measure (PRM) on $\RR_+\times \RR^n$, with intensity measure $ds\, \lambda(dz)$,   $\lambda$ is a $\sigma$-finite measure on $\RR^n$, $b, \sigma$ satisfy (D1)-(D3) as before, $\psi: \RR_+^d\times \RR^n \to \RR^d$ is a suitable coefficient map,  $h: \RR^d \to \RR^d$ is a nice truncation function (see \cite{atabudEJP}) and $h'(\x)= \x-h(\x)$ for $\x\in \RR^d$. 

The work \cite{atabudEJP}  shows that the cone $\clc$ can be used to study stability for a very general class of constrained jump-diffusions of the form in \eqref{eq:genrefdif3} for which only the first moment is given to be finite. Recall that, in the case of constrained diffusions as in \eqref{eq:genrefdif}, the condition for stability was  Assumption $\clc$-delta on the drift vector field $b(\cdot)$.
Here, in the case of a jump-diffusion, the natural definition of the drift vector field is  $\tilde \beta := 
 \cll \mbox{id}$, where $\cll$ denotes the generator of the “unconstrained” jump-diffusion away from the 
 boundary, and $\mbox{id}$ denotes the identity mapping on $\RR^d$. In the special case of a L\'{e}vy 
 process with finite mean, the drift is simply $\tilde \beta(\x) = E_{\x}(X(1))-\x$ (which is independent of 
 $\x$) where $E_{\x}$ denotes the expectation under which $X$ starts from $\x$. The basic stability assumption in  \cite{atabudEJP} is that the range of $\tilde \beta$ is 
 contained in $\cup_{k\in \NN} k \clc_1$ where $\clc_1$ is a compact subset of the interior of $\clc$. 
 Under this assumption it is shown in \cite{atabudEJP} that
 there exists a compact set $A$ such that for any compact $C \subset \RR_+^d$,
$\sup_{\x\in C} E_{\x}(\tau_A)<\infty$,  where
$\tau_A$ is the first time $X$ hits $A$. The proof of this result is based on the construction of a Lyapunov function, and on a detailed separate analysis of small and large jumps of the Markov process. As another consequence of the existence of a Lyapunov function it is shown in  \cite{atabudEJP} that $X$ is bounded in probability. From the Feller property of the process it then follows that it admits at least one invariant measure. Finally, under further additional communicability conditions  it follows that the Markov process is positive Harris recurrent and admits a unique invariant measure (see \cite{atabudEJP}).

Up to now, stability results discussed in this section were for a setting where Property \ref{prop:regsp} 
is satisfied. In the setting of a $(\mu, \D, R)$-SRBM, under the completely- $\mathcal{S}$ condition on $R$, the ergodicity 
behavior has been studied in \cite{DW}. In this work the main stability condition is formulated in terms 
of constrained versions of the linear trajectory $\eta(t;\x) = \x + \mu t$, $t\ge 0$, $\x\in \RR_+^d$. Since the 
associated Skorohod problem is not necessarily well-posed in the setting of \cite{DW} there may be multiple 
constrained versions of $\eta(t;\x)$ that are consistent with the reflection mechanism specified by $R$. 
The paper \cite{DW} shows that if all these constrained paths $\theta(t;\x)$, for all initial conditions 
$\x$ are attracted to the origin (namely $\theta(t;\x) \to 0$ as $t\to \infty$) then the SRBM is positive 
recurrent and has a unique stationary distribution.  The basic idea in the proof is to construct a 
suitable $C^2$-Lyapunov function $W$. A key property of this Lyapunov function,  in addition to
certain conditions on directional derivatives of $W$ at the boundary points, is that, for some $c>0$, 
$\nabla W(\x) \cdot \mu \le -c <0$ for all $\x \in \RR_+^d\setminus \{0\}$ . This property allows one to construct 
suitable supermartingales from which uniform in time moment estimates can be obtained.
In fact this Lyapunov function was also the inspiration for the Lyapunov function, for the general state dependent drift and a jump-diffusion setting, constructed in \cite{atabudEJP}. There, due to the state dependence,  the Lyapunov function $\tilde W$ needs to satisfy the property $\nabla \tilde W \cdot v \le -c<0$, uniformly,  for all $v$ in the range of the drift $\tilde \beta =
 \cll \mbox{id}$ which makes the construction  more involved.

\subsection{Geometric Ergodicity.}
Above results focused on existence and uniqueness of stationary distributions, and positive recurrence. A key question of interest is how long it takes for the RBM started from arbitrary initial distributions to approach stationarity (rate of convergence).  These questions were studied in \cite{budlee}. This work studied $(\mu, \D, R)$-SRBM models where $R$ is completely-$\cls$ that satisfies the stability condition of \cite{DW}, and constrained diffusion of the form in \eqref{eq:genrefdif2} with coefficients satisfying (D1)-(D3), the associated reflection matrix $R$ satisfying Property \ref{prop:regsp}, and drift satisfying the stability condition ($\clc$-delta). 
For both of these settings \cite{budlee} constructs a suitable exponentially growing Lyapunov function $V$ and establishes that the processes are $V$-uniformly ergodic. The property of $V$-uniform ergodicity says in particular that law of the process at any fixed time $t$ converges to the stationary distribution, at an exponential rate, in the total variation distance as $t \rightarrow \infty$. In fact since $V$ is exponentially growing, from this convergence one can also deduce convergence of integrals of unbounded functions (with possible exponential growth) with respect to the law at time $t$, to the corresponding integrals with respect to the stationary distribution, at an exponential rate, as $t \to \infty$.
 The result of $V$-uniform ergodicity is also used in \cite{budlee} to prove that the unique invariant measures of the SRBM and of the reflected diffusions in  \eqref{eq:genrefdif2} admit  finite moment generating functions in suitable neighborhoods of zero. As other consequences of this ergodicity one can obtain uniform (in time and initial condition in a compact set) estimates on exponential moments of these constrained processes. One can also give growth estimates of polynomial moments of the process as a function of the initial condition and
 obtain functional central limit theorems 
 for processes $\xi_{n}(t) := \frac{1}{\sqrt{n}}\left(\int_{0}^{n t}\left[F\left(X(s)\right)-\int F d\pi)\right] \mathrm{d} s\right),$ as $n\to \infty$, where $\pi$ is the unique invariant measure and $F$ is allowed to have exponential growth, and characterize the asymptotic variance  via the solution of a related Poisson equation (see\cite{budlee} for these results).
 The key steps in the proofs are the construction of a suitable Lyapunov function and establishing a minorization condition on a sufficiently large compact set (referred to as a `small set'). The Lyapunov function provides good control on the exponential moments of the return times to the small set while the minorization condition implies the existence of abstract couplings of two copies of the process (via construction of `pseudo-atoms' as described in Chapter 5 of \cite{meyn2012markov}) which have a positive chance of coalescing inside the small set. Together, they furnish exponential rates of convergence (in a weighted total variation distance).

Other than \cite{budlee}, only a few works obtain quantitative rates of convergence in some special cases. The paper \cite{IPS} obtains explicit convergence rates for time averages of bounded functionals of the state process for a class of reversible rank-based diffusions with explicit stationary measures using Dirichlet form techniques. 
The setting of one-dimensional RBM is considered in \cite{wang2014measuring} where (among other results) an estimate on the spectral gap is provided as a function of the drift and the diffusion coefficient.

\section{Parameter and 
dimension dependence on rates of convergence}\label{con}

The results of \cite{budlee}, due to the somewhat implicit treatment of the process inside the small set, provide convergence  rates that shed little light on how they qualitatively  depend on the system parameters or the state dimension. 
In this section and the next the focus is on studying, how does the convergence rate to stationarity depend on the initial distribution, diffusion parameters and the underlying dimension $d$? Moreover, it has been empirically observed that for certain types of RBM with parameters $(\mu^{(d)},\Sigma^{(d)}, R^{(d)})$ indexed by the dimension $d\ge 1$, for any fixed $k \in \mathbb{N}$, the marginal distribution of the first $k$-coordinates approaches the stationary marginal at a \emph{rate independent of $d$}. We will present results that identify conditions on $(\mu^{(d)}, \Sigma^{(d)}, R^{(d)})$ that give rise to this phenomenon, referred to as `dimension-free' convergence, and give quantitative bounds on this rate. These results set the stage for studying ergodicity properties of \emph{infinite-dimensional reflected diffusions} which will be the focus of the last two sections. Among other things, a key challenge in studying long-time behavior of such diffusions lies in the fact that they can have multiple  stationary distributions. Thus, even before making sense of rates of convergence, one needs to understand which stationary distribution the process will approach starting from a given initial distribution.

In a recent work, \cite{BC2016} obtains dimension dependent bounds on rates of convergence in Wasserstein distance for a certain class of RBM whose parameters satisfy certain `uniformity conditions' in dimension (see (BC1)-(BC3) in Section \ref{sec:bcrbm}). Although the assumptions imposed are somewhat rigid, there is a key new idea there which is to replace the role of abstract couplings (that for example underlie the approach in \cite{budlee}) with more tractable \emph{synchronous couplings} (namely, couplings where the RBM starting from different points are driven by the same Brownian motion) to analyze the behavior of the RBM inside the small set. Using explicit couplings to obtain better convergence rate estimates is a relatively recent but developing area. See \cite{bolley2010trend,Eberle2015,Eberle2017,eberle2016quantitative} for such results for other classes of diffusions. In this Section, we present results from  \cite{BB2020} which, building on the work of \cite{BC2016}, systematically develop the application of synchronous couplings and construct associated Lyapunov functions to obtain quantitative rates of convergence for any RBM satisfying Assumptions (A1)-(A3). These rates explicitly highlight the dependence of system parameters and dimension. Further, the obtained rates are significantly better when applied to the class of RBM considered in \cite{BC2016} (see Section \ref{sec:bcrbm}).

\subsection{Main Result}

We now present the main result of \cite{BB2020}. Given probability measures $\mu$ and $\nu$ on $\mathbb{R}^d_+$, a probability measure $\gamma$ on $\mathbb{R}^d_+\times \mathbb{R}^d_+$
is said to be a coupling of $\mu$ and $\nu$ if $\gamma(\cdot \times \mathbb{R}^d_+) = \mu(\cdot)$ and $\gamma(\mathbb{R}^d_+ \times \cdot) = \nu(\cdot)$.
 The $L^1$-Wasserstein distance between two probability measures $\mu$ and $\nu$ on $\mathbb{R}^d_+$ is given by
$$
W_1(\mu, \nu) = \inf\left\lbrace\int_{\mathbb{R}^d_+ \times \mathbb{R}^d_+}\|\x-\y\|_1 \gamma(\operatorname{d}\x, \operatorname{d}\y) \ : \gamma \text{ is a coupling of } \mu \text{ and } \nu\right\rbrace,
$$
where for a vector $z \in \mathbb{R}^d$, $\|z\|_1 = \sum_{i=1}^d |z_i|$.
We will denote the law of a random variable $X$ by $\mathcal{L}(X)$. 
Recall that from \cite{harrison1987brownian}, under Assumptions (A1)-(A3), there is a unique stationary distribution
of the RBM.
Denote by $\mathbf{X}(\infty)$ a random vector sampled from this stationary distribution.
Define the \emph{relaxation time}, $t_{rel}(\x)$  for the RBM starting from $\x \in \mathbb{R}^d_+$ as 
$$t_{rel}(\x) := \inf\{t\ge 0: W_1\left(\mathcal{L}(X(t; \x)), \mathcal{L}(\mathbf{X}(\infty))\right) \le 1/2\}.$$
We will abbreviate the parameters of the RBM as $\Theta :=(\mu, \Sigma, R)$. Recall that these parameters are required to satisfy (A1)-(A3). We will quantify rate of convergence to equilibrium in terms of the following functions of $\Theta,d$. 
Define the {\em contraction coefficient}
\begin{equation}\label{nd}
n(R) := \inf\{ n \ge 1: \|P^n \mathbf{1}\|_{\infty} \le 1/2\},
\end{equation}
where $\mathbf{1}$ is a $d$-dimensional vector of ones and for $u \in \mathbb{R}^d$, $\|u\|_{\infty} := \sup_{1\le i \le d} |u_i|$.
By Assumption (A1), $n(R) < \infty$.

Fix $\kappa \in (0,\infty)$. For any $\x \in \mathbb{R}^d_+$, define $\|\x\|_{\infty}^* := \sup_{1\le i \le d}\sigma_i^{-1}\x_i$. Let
\begin{align*}
a(\Theta) &:= \sup_{1 \le i \le d}\left[\frac{\sum_{j=1}^d(R^{-1})_{ij}\sigma_j}{b_i}\right],\ \ \ \ \ \ \ 
b(\Theta) := \sup_{1\le i \le d}\left[\frac{\sum_{j=1}^d(R^{-1})_{ij}\sigma_j}{\sigma_i}\right],\\
R_1(\Theta,d) &:= n(R)(1+a(\Theta)^2\log(2d)),\ \ \ \ \ \ \ 
R_2(\Theta) := a(\Theta)^2b(\Theta),\\
C_1(\x,\Theta)&:= 2\|\x\|_1 + a(\Theta)\sum_{i,j}(R^{-1})_{ij}\sigma_j,\\
C_2(\x,\Theta, \kappa) &:= 2\|\x\|_1e^{3(\kappa a(\Theta)b(\Theta))^{-1}\|\x\|_{\infty}^*} + a(\Theta)\left[2d(1+d)\left(\sum_{i,j}(R^{-1})^2_{ij}\right)\left(\sum_{j=1}^d\sigma_j^2\right)\right]^{1/2}.
\end{align*}
\begin{theorem}[\cite{BB2020}]\label{wasthm}
There exist a  $t_0 \in (0,\infty)$ and $D_1, D_2 \in (0,\infty)$ such that for every $d \in \mathbb{N}$, 
$\x \in \mathbb{R}^d_+$, every parameter choice $\Theta$, and $ t \ge t_0 \left(1 + (a(\Theta))^2\log(2d)\right)$,
\begin{align*}
	W_1\left(\mathcal{L}(X(t; \x), \mathcal{L}(\mathbf{X}(\infty))\right) &\le
\mathbb{E}(\|X(t;\x) - X(t; \mathbf{X}(\infty))\|_1)\\
& \le C_1(\x,\Theta)\left(2e^{-\frac{D_1t}{R_1(\Theta,d)}} + e^{-\frac{t}{16D_2R_2(\Theta)}}\right)
 + C_2(\x,\Theta, D_2)e^{-\frac{t}{8D_2R_2(\Theta)}}.
\end{align*}
 In particular, the relaxation time satisfies
\begin{align*}
t_{rel}(\x) \le \max\{D_1^{-1}R_1(\Theta,d) \log(8C_1(\x,\Theta)) + 16D_2R_2(\Theta)\log[4(C_1(\x,\Theta) + C_2(\x,\Theta, D_2))],\\
\qquad \qquad \qquad t_0 \left(1 + (a(\Theta))^2\log(2d)\right)\}.
\end{align*}
\end{theorem}

\subsection{Outline of Approach}
We now give a very broad outline of the approach taken in \cite{BB2020} for proving Theorem \ref{wasthm}.
\begin{itemize}
\item[(i)] As mentioned before, the key idea is to use synchronous couplings, namely, coupled RBMs from different initial states but driven by the same Brownian motion. An important feature of such couplings is that the $L^1$-distance between the two processes $X(\cdot;\0)$ and $X(\cdot;\x)$ is non-increasing in time. The rate of decay of this $L^1$-distance is quantified in terms of the contraction coefficient defined in \eqref{nd}. Roughly speaking,  the process is observed at random times $\{\eta^k(\x) : k \ge 1\}$ such that all the co-ordinates of $X(\cdot;\x)$ have hit zero at least once between successive $\eta^k(\x)$'s. It turns out that the $L^1$-distance between the synchronously coupled processes decreases by a factor of $1/2$ by time $\eta^{n(R)}(\x)$. Hence, to obtain bounds on $L^1$-Wasserstein distance, it suffices to estimate how many of these $\eta^k(\x)$'s are smaller than some given large time $t$.
\item[(ii)] For any $\mathbf{v}>0$ in $\mathbb{R}^d$ satisfying $R^{-1}\mathbf{v} \le \bb$, the RBM with oblique reflection in \eqref{RBMdef} can be dominated in a suitable manner by a normally reflected Brownian motion with drift $-\mathbf{v}$ and the same driving Brownian motion. One can construct appropriate Lyapunov functions and associated small sets for this dominating process by analyzing weighted sup-norms of its co-ordinates. This can be then used to estimate the number of times the co-ordinates of the dominating process hit zero by time $t$ as a function of $\mathbf{v}$. This, in turn, lower bounds the number of $\eta^k(\x)$'s for the original RBM less than $t$ in terms of $\mathbf{v}$. Optimizing over feasible choices of $\mathbf{v}$ gives the result.
\end{itemize}

\subsection{Examples}\label{eg}
We now describe how Theorem \ref{wasthm} can be used to obtain bounds on the rate of convergence to equilibrium in two examples that are discussed in Sections \ref{sec:bcrbm} and \ref{sec:atlmod} below.

\subsubsection{\textbf{Blanchet-Chen RBM}}
\label{sec:bcrbm}
This refers to the class of RBM under the set of assumptions in \cite{BC2016}, which are:
\begin{itemize}
\item[(BC1)] The matrix $P$ is substochastic and there exist $\kappa>0$ and $\beta \in (0,1)$ not depending on the dimension $d$ such that $\|\one^TP^n\|_{\infty} \le \kappa(1 - \beta)^n$ for all $n \ge 0$. 
\item[(BC2)] There exists $\delta>0$ independent of $d$ such that $R^{-1}\mu < -\delta \one$.
\item[(BC3)] There exists $\sigma>0$ independent of $d$ such that $ \sigma_i :=\sqrt{\Sigma_{ii}}$ satisfies $\sigma^{-1} \le \sigma_i \le \sigma$ for every $1 \le i \le d$.
\end{itemize}
Under the above conditions \cite{BC2016} give a polynomial bound  of $O(d^4(\log d)^2)$  on the relaxation time of the RBM.
As shown in the following theorem, Theorem \ref{wasthm} gives a  substantial improvement by establishing a polylogarithmic relaxation time of $O((\log d)^2)$. 
\begin{theorem}[\cite{BB2020}]\label{BCbound}
Under Assumptions (BC1), (BC2) and (BC3), there exist positive constants $E_1, E_2, E_3, E_4,t_1$ such that for any $\x \in \mathbb{R}^d_+, t \ge t_1 \max\{\|\x\|_{\infty}, \log(2d)\}$,
\begin{equation*}
\mathbb{E}(\|X(t;\x) - X(t; \mathbf{X}(\infty))\|_1) \le  2\left(2\|\x\|_1 + E_1d^2\right)e^{-E_2 t/\log(2d)} + \left(4\|\x\|_1 + E_1d^2\right)e^{-E_4t/2} + E_3 d^{2}e^{-E_4t}.
\end{equation*}
In particular, the relaxation time satisfies
\begin{align*}
t_{rel}(\x) \le \max \left\lbrace E_2^{-1}\log\left[8\left(2\|\x\|_1 + E_1d^2\right)\right]\log (2d) + E_4^{-1}\left[2\log\left[8\left(4\|\x\|_1 + E_1d^2\right)\right] + \log (8E_3d^{2})\right],\right.\\
\left.\qquad \qquad \qquad  t_1 \max\{\|\x\|_{\infty}, \log(2d)\}\right\rbrace.
\end{align*}
\end{theorem}

\subsubsection{\textbf{Gap process of rank-based diffusions}} 
\label{sec:atlmod}
Rank-based diffusions are interacting particle systems where the drift and diffusion coefficient of each particle depends on its rank. Mathematically, they are represented by the SDE:
\begin{equation}\label{rbdef}
dY_i(t) = \left(\sum_{j=0}^{d} \delta_j \one_{[Y_i(t) = Y_{(j)}(t)]}\right)dt + \left(\sum_{j=0}^{d} \sigma_j \one_{[Y_i(t) = Y_{(j)}(t)]}\right)dW_i(t)
\end{equation}
for $0\le i \le d$, where $\{Y_{(j)}(t) : t \ge 0\}$ denotes the trajectory of the rank $j$ particle as a function of time $t$ ($Y_{(0)}(t) \le \dots\le Y_{(d)}(t) \text{ for all } t \ge 0$), $\delta_j, \sigma_j$ denote the drift and diffusion coefficients of the rank $j$ particle, and
$W_i$, $0\le i \le d$, are mutually independent standard one dimensional Brownian motions.
We will assume throughout that $\sigma_i>0$ for all $0 \le i \le d$. Rank-based diffusions have been proposed and extensively studied as models for problems in finance and economics. A special case is the  Atlas model \cite{Fern}  where the minimum particle (i.e. the particle with rank $1$) is a Brownian motion with positive drift and the remaining particles are Brownian motions without drift (i.e. $\delta_i=0$ for all $i\ge1$). The general setting considered in 
\eqref{rbdef} was introduced in 
\cite{BFK}. 
In order to study the long time behavior, it is convenient to consider the gap process $Z = (Z_1,\dots,Z_{d})$, given by $Z_i = Y_{(i)} - Y_{(i-1)}$ for $1 \le i \le d$. 
The process $Z \equiv Z(t;\z)$ is a RBM in $\mathbb{R}_+^{d}$ given as
\begin{equation*}
Z(t;\z) = \z + \D B(t) + \mu t +  RL(t)
\end{equation*}
where $\z$ is the initial gap sequence, $B$ is a standard $d$-dimensional Brownian motion, $\mu_i = \delta_{i} - \delta_{i-1}$ for $1 \le i \le d$, $\D \in \mathbb{R}^{d \times d}$, $L$ is the local time process 
associated with $Z$ and $R$ satisfies Assumption (A1), namely it belongs to the Harrison-Reiman class. The covariance matrix $\Sigma = \D\D^T$ has entries $\Sigma_{ii} = \sigma_{i-1}^2 + \sigma_{i}^2$ for $1 \le i \le d$, $\Sigma_{i(i-1)} = - \sigma_{i-1}^2$ for $2 \le i \le d$, $\Sigma_{i(i+1)} = -\sigma_{i}^2$ for $1 \le i \le d-1$ and $\Sigma_{ij} = 0$ otherwise. In particular, (A3) is satisfied, namely $\Sigma$ is positive definite.
Moreover, $R$ is given explicitly as 
$R = I-P^T$, where $P$ is the substochastic matrix given by $P_{i (i+1)} = P_{i (i-1)} = 1/2$ for all $2 \le i \le d-1$, $P_{12} = P_{d(d-1)} = 1/2$ and $P_{ij} = 0$ if $|i-j| \ge 2$. 
 From \cite{harrison1987brownian} the process is positive recurrent and has a unique stationary distribution
if 
Assumption (A2) is satisfied, namely $\bb =-R^{-1}\mu>0$, which  is  same as the following condition:
\begin{equation}\label{atlasstab}
b_k = \sum_{i=1}^k(\delta_{i-1} - \overline{\delta})>0 \mbox{ for } 1 \le k \le d, \mbox{ where } \overline{\delta} = (d+1)^{-1}\sum_{j=0}^{d} \delta_j.
\end{equation}
In the special case where
\begin{equation}\label{sk}
\sigma_{i}^2 - \sigma_{i-1}^2 = \sigma_{1}^2 - \sigma_0^2 \text{ for all } 1 \le i \le d,
\end{equation}
 the stationary distribution is explicit and
 takes the form $\mathcal{L}(\mathbf{Z}(\infty)) = \otimes_{k=1}^d \operatorname{Exp}(2b_k\left(\sigma_{k-1}^2 + \sigma_{k}^2\right)^{-1})$ (see Section 5 of \cite{ichiba2011hybrid}). 
 For the general case (i.e. $\sigma_i$ are strictly positive and \eqref{atlasstab} is satisfied) explicit formulas for stationary distribution are not available, however 
 from \cite{{budlee}}, the law of $Z(t;\z)$ converges to the unique stationary distribution in (weighted) total variation distance at an (albeit implicit) exponential rate.  

From Theorem \ref{wasthm} one can give the following bound on the rate of $L^1$-Wasserstein convergence of the gap process to $\mathbf{Z}(\infty)$.
Note that we do not  require reversibility or an explicit expression for the stationary measure. 

Two key quantities appearing in the rate of convergence are
\begin{equation}\label{eq:sigupplow}
a^* := \sup_{1 \le i \le d}\frac{i(d+1-i)}{b_i}, \ \ \ \sigma = \left(\sup_{0 \le i \le d} \sigma_i \right) \vee \left(\sup_{0 \le i \le d} \sigma_i^{-1} \right)
\end{equation}
where $b_i$ are defined in \eqref{atlasstab} and $\sigma_i$ is the standard deviation of the rank $i$ particle (see \eqref{rbdef}).
\begin{theorem}[\cite{BB2020}]\label{atlasshrugged}
There exist positive constants $F_1, F_2, F_3, F_4, t_2$ such that for any $d\in \NN$, $\z \in \mathbb{R}^d_+$ and any $t \ge t_2 \max\{\sigma^2a^*\|\z\|_{\infty}, 1 + \sigma^2a^{*2}\log(2d)\}$,
\begin{align*}
\mathbb{E}(\|Z(t;\z) - Z(t; \mathbf{Z}(\infty))\|_1) \le 2\left(2\|\z\|_1 + F_1\sigma^2a^*d^3\right)e^{-F_2 t/\left[d^2(1+\sigma^2a^{*2}\log(2d))\right]}\\
 + \left(4\|\z\|_1 + F_1\sigma^2a^*d^3\right)e^{-F_4t/[2\sigma^4a^{*2}(d+1)^2]} + F_3\sigma^2a^*d^{7/2}e^{-F_4t/[\sigma^4a^{*2}(d+1)^2]}.
\end{align*}
In particular, the relaxation time satisfies
\begin{align*}
t_{rel}(\z) &\le \max \left\lbrace F_2^{-1}\left[d^2(1+ \sigma^2a^{*2}\log(2d))\right]\log\left[8\left(2\|\z\|_1 + F_1\sigma^2a^*d^3\right)\right]\right.\\
&\left. \qquad \qquad\qquad + F_4^{-1}\sigma^4a^{*2}(d+1)^2\left[2\log\left[8\left(4\|\z\|_1 + F_1\sigma^2a^*d^3\right)\right] + \log (8F_3\sigma^2a^*d^{7/2})\right],\right.\\
&\left. \qquad \qquad\qquad\qquad t_2 \max\{\sigma^2a^*\|\z\|_{\infty}, 1 + \sigma^2a^{*2}\log(2d)\}\right\rbrace.
\end{align*}
\end{theorem}
\begin{rem}\label{sam}
The standard Atlas model \cite{Fern} is a special case of \eqref{rbdef} with $\delta_0=1$, $\delta_i = 0$ for all $i \ge 1$ and $\sigma_i=1$ for all $i$. For this model, using \eqref{atlasstab}, for any $k \ge 1$,
$$
b_k = \sum_{i=1}^k(\delta_{i-1} - \overline{\delta}) = \frac{(d+1-k)}{d+1}
$$
and
$$
a^* := \sup_{1 \le i \le d}\frac{i(d+1-i)}{b_i} = \sup_{1 \le i \le d}i(d+1) = d(d+1), \ \ \ \sigma=1.
$$
Using these in Theorem \ref{atlasshrugged}, one obtains positive constants $G_1, G_2, G_3, G_4, t_3$ such that for any $d\in \NN$, $\z \in \mathbb{R}^d_+$ and any $t \ge t_3\{d^2\|\z\|_{\infty}, 1 + d^2\log(2d)\}$,
$$
\mathbb{E}(\|Z(t;\z) - Z(t; \mathbf{Z}(\infty))\|_1) \le G_1\left(\|\z\|_1 + d^5\right)e^{-G_2t/d^6\log(2d)} + G_3d^{11/2}e^{-G_4t/d^6}.
$$
In particular, the relaxation time for the standard Atlas model is $O(d^6(\log d)^2)$ as $d \rightarrow \infty$.
\end{rem}

\section{Dimension-free convergence rates for local functionals of RBM}
Typically, growing dimension slows down the rate of convergence for the whole system, as is reflected in the bounds obtained in the previous section, but empirically one often observes a much faster convergence rate to equilibrium of \emph{local statistics} of the system. 
In this Section, we describe results for a class of RBMs for which convergence rates of local statistics do not depend on the underlying dimension of the entire system. More precisely, consider a family of processes $X^{(d)} \sim \RBM(\mu^{(d)}, \Sigma^{(d)}, R^{(d)})$ indexed by the dimension $d\ge 1$ (the superscript $(d)$ is dropped subsequently for notational convenience). Conditions on $(\mu^{(d)}, \Sigma^{(d)}, R^{(d)})$ are sought under which, for any fixed $k \in \mathbb{N}$, the marginal distribution of the first $k$-coordinates approaches the stationary marginal at a rate independent of $d$. The  results below are taken from \cite{banerjee2020dimension}, where this phenomenon is named \emph{dimension-free local convergence}.

Mathematically, this is challenging as the local evolution is no longer Markovian and the techniques in \cite{BC2016,BB2020} cannot be readily applied. A crucial observation is that certain \emph{weighted $L^1$-distances} (see $\|\cdot\|_{1,\beta}$ defined in Section \ref{wnormdef}) between synchronously coupled RBMs show dimension-free contraction rates. The evolution of such weighted distances are tracked in time for synchronously coupled RBMs $X(\cdot; \0)$ and $X(\cdot;\x)$ for $\x \in \mathbb{R}^d_+$. It is shown in \cite{banerjee2020dimension} that for this distance to decrease by a dimension-free factor of its original value, only a subset of co-ordinates of $X(\cdot;\x)$, whose cardinality depends on the value of the original distance (and not directly on the dimension $d$), need to hit zero. This is in contrast with the unweighted $L^1$-distance considered in \cite{BC2016,BB2020} where all the coordinates need to hit zero to achieve such a contraction, thereby slowing down the convergence rate. Consequently, by tracking the hitting times to zero of a time dependent number of co-ordinates, one achieves dimension-free convergence rates in this weighted $L^1$-distance as stated in Theorem \ref{thm:main_fromx}. This, in turn, gives dimension-free local convergence as is made precise in \eqref{locstat}. In Section \ref{sec:examples}, an application of Theorem \ref{thm:main_fromx} is presented for the gap process of a rank-based interacting particle system with asymmetric collisions to obtain explicit rates exhibiting dimension-free local convergence. 

\subsection{A weighted $L^1$-distance governing dimension-free local convergence}\label{wnormdef}
The investigation of dimension-free convergence relies on the analysis of the weighted $L^1$-distance $\|X(\cdot;\x) - X(\cdot;\XX(\infty))\|_{1,\beta} := \sum_{i=1}^d \beta^i |X_i(\cdot; \x) - X_i(\cdot; \XX(\infty))|$ in time, for appropriate choices of $\beta \in (0,1)$.
Towards this end, the following functionals play a key role:
\begin{align}
    u_{\beta}(\x, t) &= \bnorm{\rinv\left(X(t; \x) - X(t; \0)\right)} := \sum_{i=1}^d \beta^i \left| \left[\rinv\left(X(t; \x) - X(t; \0)\right)\right]_i \right|,\label{eqn:weighted_norm}\\ 
    u_{\pi,\beta}(t) &= u_{\beta}(\XX(\infty), t), \quad \quad t \ge 0. \label{eqn:u_pi}
\end{align}
In the following, when $\beta$ is clear from context, we will suppress dependence on $\beta$ and write $u$ for $u_{\beta}$ and $u_{\pi}$ for $u_{\pi,\beta}$. The above functionals are convenient because the vector $\rinv\left(X(t; \x) - X(t; \0)\right)$ is co-ordinate wise non-negative and non-increasing in time. Moreover, $\rinv\left(X(t; \x) - X(t; \0)\right) \ge X(t; \x) - X(t; \0) \ge 0$ for all $t \ge 0$. This fact and the triangle inequality can be used to show for any $\x \in \mathbb{R}^d_+, t \ge 0$,
\begin{align*}
\bnorm{\left(X(t; \x) - X(t; \XX(\infty))\right)} \le u(\x, t) + u_\pi(t),
\end{align*}
where the bound on the right hand side is non-increasing in time. Due to the monotonicity of the bound, it suffices to find `events' along the trajectory of the coupled processes that lead to a reduction in this bound by a dimension-independent factor. Using this idea, conditions are obtained under which there exists a $d$-independent $\beta \in (0,1)$ and a function $f: \Reals_+ \mapsto \Reals_+$ not depending on the dimension $d$ such that $f(t) \to 0$ as $t \to \infty$ and, for any $\x$ in an appropriate subset $\mathcal{S}$ of $\mathbb{R}_+^d$,
\begin{equation}\label{eqn:dimension_free}
    \Expect{\bnorm{\left(X(t; \x) - X(t; \XX(\infty))\right)}} \le  \Expect{u(\x, t) + u_\pi(t)} \le C f(t), \quad \quad t \ge t_0,
\end{equation}
where $C, t_0 \in (0, \infty)$ are constants not depending on $d$ (but can depend on $\x$). This, in particular, gives dimension-free local convergence in the following sense: For any $k \in \{1,\dots,d\}$, consider any function $\phi: \mathbb{R}^k_+ \mapsto [0,\infty)$ which is Lipschitz, i.e., there exists $L_{\phi} >0$ such that
$$
|\phi(\mathbf{u}) - \phi(\mathbf{v})| \le L_{\phi} \|\mathbf{u}-\mathbf{v}\|_1, \quad \quad \mathbf{u}, \mathbf{v} \in \mathbb{R}^k_+.
$$
Recall that $\mathcal{L}(Z)$ denotes the law of a random variable $Z$ and $W_1$ denotes the associated $L^1$-Wasserstein distance. Then, \eqref{eqn:dimension_free} implies for $\x \in \mathcal{S}$,
\begin{align}\label{locstat}
W_1\left(\mathcal{L}(\phi(X|_k(t; \x))), \mathcal{L}(\phi(\XX|_k(\infty)))\right) 
&\le \Expect{|\phi(X|_k(t; \x)) - \phi(X|_k(t; \XX(\infty)))|} \nonumber \\
&\le C\beta^{-k} L_{\phi} f(t), \quad \quad t \ge t_0,
\end{align}
where for $\x = (x_1, \ldots , x_d)' \in \RR_+^d$, $\x|_k = (x_1, \ldots, x_k)$.
\subsection{Parameters and Assumptions}
We now define the parameters that govern dimension-free local convergence which, in turn, are defined in terms of the original model parameters $(\mu,\Sigma,R)$ of the associated RBM. 
We abbreviate $\sigma_i = \sqrt{\Sigma_{ii}}, i=1,\ldots, d$. Define for $1 \le k \le d$,
\begin{align}\label{pardefcon}
\bk{k} &:= -\rinvk{k}\muk{k}, \quad \quad \bb = \bk{d}, \notag\\
\blowk{k} &:= \min_{1 \le i \le k} \bk{k}_i,  \quad \quad \ak{k} := \max_{1\le i \le k} \frac{1}{\bk{k}_i} \sum_{j=1}^k (\rinvk{k})_{ij}\sigma_j,
\end{align}
where $\muk{k}$ denotes the restriction of $\mu$ onto first $k$ co-ordinates and $\rinvk{k}$ denotes the inverse of the principal $k \times k$ submatrix of $R$. 
To get a sense of why these parameters are crucial, recall that our underlying strategy is to obtain contraction rates of $u(\x,\cdot)$ defined in \eqref{eqn:weighted_norm} by estimating the number of times a subset of the co-ordinates of $X(\cdot;\x)$, say $\{X_1(\cdot; \x),\dots,X_k(\cdot; \x)\}, k \le d$, hit zero. However, this subset does not evolve in a Markovian way. Thus, we use monotonicity properties of RBMs to couple this subset with a $\mathbb{R}^k_+$-valued reflected Brownian motion $\bar{X}(\x|_k,\cdot)$, started from $\x|_k$ and defined in terms of $\mu|_k, \Sigma|_k, R|_k$ and (a possible restriction of) the same Brownian motion driving $X(\cdot;\x)$, such that $X_i(\x,t) \le \bar{X}_i(\x|_k,t)$ for all $1 \le i \le k$. The analysis in \cite{BB2020} shows that the parameters defined in \eqref{pardefcon} with $k=d$ can be used to precisely estimate the minimum number of times all co-ordinates of $X(\cdot;\x)$ hit zero by time $t$ as $t$ grows. Applying this approach, for any $1 \le k \le d$, the parameters \eqref{pardefcon} can be used to quantify analogous hitting times for the process $\bar{X}(\cdot; \x|_k)$ which, by the above coupling, gives control over corresponding hitting times of $\{X_1(\cdot; \x),\dots,X_k(\cdot; \x)\}$.

Given below is a set of assumptions on the model parameters $(\mu,\Sigma,R)$ which guarantee dimension-free local convergence.\\

\textbf{Assumptions DF:}
There exist $d$-independent constants $\sigmalow, \sigmahigh, b_0 > 0$, $r^* \ge 0$, $M, C \ge 1$, $k_0 \in \{2, \ldots, d\}$ and $\alpha\in (0, 1)$ such that for all $d \ge k_0$, 
\begin{itemize}
    \item[I.] $(\rinv)_{ij} \le C \alpha^{j-i}$ for $1\le i \le j \le d$,
    \item[II.] $(\rinv)_{ij} \le M$ for $1 \le i, j \le d$,
    \item[III.] $\blowk{k} \ge b_0 k^{-r^*}$ for $k = k_0, \ldots, d$,
    \item[IV.] $\sigma_i \in [\sigmalow, \sigmahigh]$ for $1 \le i \le d$.
\end{itemize}
We explain why Assumptions DF are `natural' in obtaining dimension-free local convergence. Since $P$ is a transient and substochastic, it can be associated to a killed Markov chain on $\{0\} \cup \{1,\dots,d\}$ with transition matrix $P$ on $\{1,\dots,d\}$ and killed upon hitting $0$ (i.e. probability of going from state $k \in \{1,\dots,d\}$ to $0$ is $1 - \sum_{l=1}^d P_{kl}$ and $P_{00} =1$). Moreover, since $P$ is transient and $R = I - P^T$, we have $\rinv = \sum_{n=0}^\infty (P^T)^n$. This representation shows that $(\rinv)_{ij}$ is the expected number of visits to site $i$ starting from $j$ of this killed Markov chain. For fixed $\x \in \mathbb{R}^d_+$ and $k < < d$, consider a local statistic of the form $\phi(X|_k(t; \x))$ as in \eqref{locstat}. For this statistic to stabilize faster than the whole system, we expect the influence of the far away co-ordinates $X|_j(\cdot; \x), j>>k,$ to diminish in an appropriate sense as $j$ increases. This influence is primarily manifested through the oblique reflection arising out of the $R$ matrix via the above representation of $R^{-1}$. I of Assumptions DF quantifies this intuition by requiring that the expected number of visits to state $i$ starting from state $j>i$ of the associated killed Markov chain decreases geometrically with $j-i$. This is the case, for example, when this Markov chain started from $j>i$ has a uniform `drift' away from $i$ towards the cemetery state. See Section \ref{sec:examples}. In more general cases, one can employ Lyapunov function type arguments \cite{meyn2012markov} to the underlying Markov chain to check I. 

II above implies that the killed Markov chain starting from state $j$ spends at most $M$ expected time at any other site $i \in \{1,\dots,d\}$ before it is absorbed in the cemetery state $0$. This expected time, as  calculations show, is intimately tied to decay rates of $\bnorm{\left(X(\cdot; \x) - X(\cdot; \XX(\infty))\right)}$. 

As noted in \cite{harrison1987brownian,BC2016,BB2020}, the `renormalized drift' vector $\bb$ characterizes positive recurrence of the whole system. Assumption III above, prescribes a power law type co-ordinate wise lower bound of the renormalized drift vector $\bb^{(k)}$ of the projected system $X|_k(\cdot; \x)$ as $k$ grows. In particular, if $\blowk{k}$ is uniformly lower bounded by $b_0$, one can take $r^*=0$.
Finally, IV above is a quantitative `uniform ellipticity' condition on the co-ordinates of the driving noise $\D B(\cdot)$.

\subsection{Main result}\label{sec:results}
The following result gives explicit bounds on the decay of expectation of the weighted distance $\|X(\cdot;\x) - X(\cdot;\XX(\infty))\|_{1,\sqrt{\alpha}}$ ($\alpha$ defined in I of Assumptions DF) with time, for RBMs satisfying Assumptions DF.
We first define some constants that will appear in Theorem \ref{thm:main_fromx}. They are needed to bound moments of certain weighted norms of the stationary random variable $\XX(\infty)$.

Suppose Assumptions DF hold, with $k_0 \in \{2, \ldots, d\}$ and $\alpha \in (0, 1)$ defined therein. Set
\begin{equation}\label{l1def}
    L_1 := k_0^{r^* + 1} +  \sum_{i = k_0}^d i^{3 + r^*} \alpha^{i/8}.
\end{equation}
Also, for $B \in (0,\infty)$, define the set
\begin{equation}\label{startset}
\mathcal{S}(b,B) := \left\lbrace x \in \mathbb{R}_+^d : \sup_{1 \le i \le d} \blow^{(i)} \supnorm{x|_i} \le B\right\rbrace.
\end{equation}
\begin{theorem}[\cite{banerjee2020dimension}]\label{thm:main_fromx}
Let $X$be a $( \mu,\Sigma, R)$-RBM.  Suppose that
Assumptions DF hold with $\alpha \in (0, 1)$ defined therein.

Fix any $B \in (0,\infty)$. Then there exist constants $C_0, C_0', C_1 > 0$ not depending on $d, r^*$ or $B$ such that with $t_0' = t_0'(r^*) = C_0'\left(1+r^*\right)^{8+4r^*}$ and $L_1, L_2, \mathcal{S}(b,B)$ as defined in \eqref{l1def}-\eqref{startset}, we have for any $\x \in \mathcal{S}(b,B)$ and any $d > t_0'^{1/(4+2r^*)}$,
\begin{multline}\label{eqn:main_fromx_result1}
    \Expect{\|X(t; \x) - X(t; \XX(\infty))\|_{1,\sqrt{\alpha}}}\\
    \le
\begin{cases}
C_1 \left(L_1\sqrt{1+t^{1/(4+2r^*)} } + \supnorm{\x}\exp\left\lbrace B/\sigmalow^2\right\rbrace\right) \> \exp\{-C_0 t^{1/(4+2r^*)} \},
  &   t_0' \le t < d^{4 + 2r^*},\\\\
C_1 \left(L_1\sqrt{1+t^{1/(4+2r^*)} } + \supnorm{\x}\exp\left\lbrace B/\sigmalow^2\right\rbrace\right)\> \exp\left\lbrace -C_0\frac{t}{d^{3 + 2r^*}}\right\rbrace,   & t \ge d^{4 + 2r^*}.
\end{cases}
\end{multline}
\end{theorem}
Theorem \ref{thm:main_fromx} directly implies dimension-free bounds on $t \mapsto \|X(t;\x) - X(t; \XX(\infty))\|_{1,\sqrt{\alpha}}$ in the sense of \eqref{eqn:dimension_free} which, in turn, produce dimension-free local convergence rates as given by \eqref{locstat}. The multi-part bound in the theorem is designed to `continuously interpolate' the \emph{dimension-free stretched exponential bounds} with the \emph{dimension-dependent exponential bounds} obtained in \cite{BB2020}.

\subsection{Application to the Asymmetric Atlas model}\label{sec:examples}
Consider the interacting particle system represented by the following SDE:
\begin{equation}\label{ascol}
Y_{(k)}(t) = Y_{(k)}(0) + \mathbf{1}_{k=1} t + B^*_k(t) + p L^*_{k}(t) - q L^*_{k+1}(t), \quad t \ge 0,
\end{equation}
for $0\le k \le d$, $p \in (0,1), q= 1-p$. Here, $L^*_{0}(\cdot) \equiv L^*_{d+1}(\cdot) \equiv 0$, and for $0 \le k \le d-1$, $L^*_{k+1}(\cdot)$ is a continuous, non-decreasing, adapted process that denotes the collision local time between the $k$-th and $(k+1)$-th co-ordinate processes of $Y_{(\cdot)}$, namely $L^*_{k+1}(0)=0$ and $L^*_{k+1}(\cdot)$ can increase only when $Y_{(k)} = Y_{(k+1)}$.
$B^*_k(\cdot)$, $0\le k \le d$, are mutually independent standard one dimensional Brownian motions. Each of the $d+1$ ranked particles with trajectories given by $(Y_{(0)}(\cdot),\dots,Y_{(d)}(\cdot))$ evolves as an independent Brownian motion (with the particle $0$ having unit positive drift) when it is away from its neighboring particles, and interacts with its neighbors through possibly asymmetric collisions. The symmetric Atlas model, namely the case $p=1/2$, was introduced in \cite{Fern} as a mathematical model for stochastic portfolio theory. The asymmetric Atlas model was introduced in \cite{karatzas_skewatlas}. It was shown that it arises as scaling limit of numerous well-known interacting particle systems involving asymmetrically colliding random walks \cite[Section 3]{karatzas_skewatlas}. Since then, this model has been extensively analyzed: see \cite{karatzas_skewatlas,IPS,ichiba2013strong,AS} and references therein.

The gaps between the particles, defined by $Z_i(\cdot) = Y_{(i)}(\cdot) - Y_{(i-1)}(\cdot), 1 \le i \le d$, evolve as an $\RBM(\mu,\Sigma, R)$ with $\Sigma$ given by $\Sigma_{ii} = 2$ for $i = 1, \ldots, d$, $\Sigma_{ij}=-1$ if $|i-j|=1$, $\Sigma_{ij} = 0$ if $|i-j| >1$, $\mu$ given by $\mu_1 = -1, \mu_j = 0$ for $j = 2, \ldots, d$, and $R = I-P^T$, where
\begin{equation}\label{eqn:skew_atlas1}
    P_{ij} = \begin{cases}
    p \quad \quad j = i+1,\\
    1-p \quad \quad j = i-1, \\
    0 \quad \quad \text{otherwise}.
    \end{cases}
\end{equation}
Clearly $R$ is in the Harrison-Reiman class. The paper \cite{banerjee2020dimension} establishes dimension-free local convergence for the gap process of the Asymmetric Atlas model when $p \in (1/2, 1)$.
Recall that the reflection matrix $R= I-P^T$ is associated with a killed Markov chain. For the Asymmetric Atlas model, this Markov chain has a more natural description as a random walk on $\{0,1,\dots,d+1\}$ which increases by one at each step with probability $p$ and decreases by one with probability $1-p$, and is killed when it hits either $0$ or $d+1$. Then for $1 \le i,j \le d$, $(\rinv)_{ij}$ is the expected number of visits to $i$ starting from $j$ by this random walk before it hits $0$ or $d+1$. Since $p > 1-p$, the random walk has a drift towards $d+1$, which suggests I, II of Assumptions DF hold. This is confirmed by direct computation, which gives for $q = 1-p$,
\begin{equation}\label{eqn:skew_atlas2}
    (\rinv)_{ij} = \begin{cases}
     \frac{\left(q/p \right)^{j-i}}{p-q}\frac{\left(1 - (q/p)^i \right)\left(1-(q/p)^{d+1-j} \right)}{1- (q/p)^{d+1}} \le \frac{\left(q/p \right)^{j-i}}{p-q} & 1 \le i \le j \le d,\\
     \frac{\left(p/q \right)^{i-j}}{p-q}\frac{\left((p/q)^j - 1 \right)\left((p/q)^{d+1-i} - 1 \right)}{(p/q)^{d+1} - 1} \le \frac{1}{p-q} & 1 \le j <  i \le d.
    \end{cases}
\end{equation}
Now I, II and IV of Assumptions DF hold with $M = C = \frac{1}{p-q}$, $\alpha = \frac{q}{p}$ and $\sigmalow=\sigmahigh=\sqrt{2}$. Furthermore, the restriction $P|_{k}$ is defined exactly as in \eqref{eqn:skew_atlas1} with $k$ in place of $d$. Thus $\rinvk{k}$ is given by \eqref{eqn:skew_atlas2} with $k$ in place of $d$, and $\bk{k} = -\rinvk{k} \mu|_k$ is the first column of $\rinvk{k}$. This entails,
\begin{align}\label{eqn:skew_atlas3}
    \bk{k}_i &= \frac{\left(p/q \right)^{i-1}}{p-q}\frac{\left((p/q) - 1 \right)\left((p/q)^{k+1-i} - 1 \right)}{(p/q)^{k+1} - 1} 
    \ge \frac{1}{q}\left(\frac{p}{q}\right)^{k-1}\frac{\left((p/q) - 1\right)}{(p/q)^{k+1}} \notag\\
    &= \frac{p-q}{p^2}  =: b_0 > 0, \quad \quad 1 \le i \le k, \> 1 \le k \le d.
\end{align}
Thus $\blowk{k} \ge b_0$ for all $1 \le k \le d$, uniformly in $d$. This shows that III of Assumptions DF holds with $b_0$ specified by \eqref{eqn:skew_atlas3} and $r^* = 0$. Moreover, it follows from the first equality in \eqref{eqn:skew_atlas3} that $\bk{k}_i  \le p/(p-q)$ for all $1 \le k \le d$ and $1 \le i \le k$. Therefore, recalling the definition of $\mathcal{S}(b,\cdot)$ from \eqref{startset}, for any $\z \in \mathbb{R}^d_+$,
$$
\z \in \mathcal{S}(b, p\supnorm{z}/(2p-1)).
$$
The above observations result in the next theorem, which follows directly from Theorem \ref{thm:main_fromx}. These bounds imply dimension-free convergence rates in the sense of \eqref{eqn:dimension_free} and \eqref{locstat}.
\begin{theorem}[\cite{banerjee2020dimension}]\label{thm:skew_atlas}
Suppose $Z$ is the RBM representing the gap process of the Asymmetric Atlas model with $p \in (1/2, 1)$. Then there exist constants $\bar{C}, \bar{C}_0, t_0' > 0$ depending on $p$ but not on $d$ such that for $d > t_0'$,
\begin{multline}
    \Expect{\|Z(t;\z) - Z(t;\ZZ(\infty))\|_{1,\sqrt{\frac{1-p}{p}}}}\\
    \le
\begin{cases}
\bar{C}\left(\sqrt{1+t^{1/4}} + \supnorm{\z} e^{p\supnorm{\z}/(4p-2)}\right)\> e^{-\bar{C}_0 t^{1/4}},
  &  t_0' \le t < d^4,\\\\
\bar{C}\left(\sqrt{1+t^{1/4}} + \supnorm{\z} e^{p\supnorm{\z}/(4p-2)}\right)\> e^{-\bar{C}_0 t/d^3},   & t \ge d^4.
\end{cases}
\end{multline}
\end{theorem}
We note  that the law of $\ZZ(\infty)$ has an explicit product form in this case. See \cite[Proposition 2.1]{AS}.

\subsection{Polynomial local convergence rates for perturbations from stationarity of the Atlas model}\label{dfsym}

A natural question is whether one can expect dimension-free local convergence when Assumptions DF do not hold. 
This is currently an open problem. For one example where one does not expect such dimension-free local convergence to hold, 
consider the Asymmetric Atlas model from Section \ref{sec:examples} when $p<1/2$.  Observe from \eqref{ascol} that a larger change of momentum happens to the lower particle (smaller $k$) during each collision. Hence, effects of unstable gaps between far-away particles can propagate in a more unabated fashion to the lower particles. Moreover, by analyzing the associated killed Markov chain, one concludes that it takes roughly linear time in $d$ for the effect of the $d$-th gap to propagate to the first gap. This heuristic reasoning suggests that one cannot expect any reasonable form of dimension-free local convergence in this case.

Now consider the `critical' case: the (symmetric) Atlas model, namely, the case $p=1/2$ in the model of Section \ref{sec:examples}. 
Note that here, $\mathbf{b} = -\rinv \mu = \{(\rinv)_{i1} \}_{i=1}^d > 0$ and $\Sigma_{ii} = 2$ for all $i$. Therefore, there exists a unique stationary distribution. In fact, if $\ZZ(\infty)$ denotes the corresponding stationary distributed gap, it holds that \cite{harrisonwilliams_expo,ichiba2011hybrid}
\begin{equation}\label{eqn:atlas_stationary}
    \ZZ(\infty) \sim \bigotimes_{i = 1}^d \text{Exp}\left(2\left(1 - \frac{i}{d+1}\right)\right).
\end{equation}
In the symmetric Atlas model, the associated killed Markov chain behaves as a simple symmetric random walk away from the cemetery points, and lacks the strong drift towards the cemetery points as formulated in I of Assumptions DF. Indeed, $\rinv$ is given as (see e.g. \cite[Proof of Theorem 4]{BB2020}, or  take $p \to 1/2$ in \eqref{eqn:skew_atlas2}):
\begin{equation}\label{eqn:atlas2}
    (\rinv)_{ij} = \begin{cases}
    2i\left(1 - \frac{j}{d+1}\right)\quad \quad 1 \le i \le j \le d,\\
    2j\left(1 - \frac{i}{d+1}\right) \quad \quad 1 \le j < i \le d.
    \end{cases}
\end{equation}
The above representation shows that $\rinv$ violates I, II of Assumptions DF (consider, for example,  $i = j = \floor{d/2}$). 
Hence, the techniques sketched above for obtaining dimension-free convergence fail to apply in this case. Instead, a different approach is taken in \cite{banerjee2020dimension} which involves analyzing the long term behavior of pathwise derivatives of the RBM in initial conditions. Using this analysis, \emph{polynomial convergence rates} to stationarity in $L^1$-Wasserstein distance are obtained when the initial distribution of the gaps between particles is in an appropriate perturbation class of the stationary measure (see \cite[Section 3.1, Definition 1]{banerjee2020dimension}). We give a representative result obtained as a special case of the more general theorem \cite[Theorem 4]{banerjee2020dimension}.

\begin{theorem}[\cite{banerjee2020dimension}]\label{cor:expo_perturbation}
Consider $\mathbf{Y} = (Y_1, Y_2, \ldots)^T$ where $\{Y_i\}_{i \ge 1}$ are independent random variables with $Y_i \sim \operatorname{Exp}(i^{1+\beta})$ (exponential with mean $i^{-(1+\beta)}$), for some $\beta > 0$. Then there exist positive constants $C_1, C_2$ not depending only on $\beta$ such that
\begin{equation*}
  \Expect{\|Z(t; \ZZ(\infty) + \mathbf{Y}\vert_d) - Z(t;\ZZ(\infty)) \|_1} \ 
    \le \
\begin{cases}
C_1t^{-\frac{\beta}{1 + \beta} \frac{3}{32}},
    & t_0' \le t < d^{16/3},\\\\
C_2d\exp \left(-C_0 \frac{t}{d^6\log (2d)} \right), 
    & t\ge t_0'' \ d^4\log (2d),
\end{cases}
\end{equation*}
where $t_0' \in (0, \infty)$ does not depend on $d$ or $\beta$.
\end{theorem}
The analysis of the derivative process is based on the analysis of a random walk in a random environment generated by the times and locations where the RBM hits faces of $\mathbb{R}^d_+$. We mention here that \cite{blanchent2020efficient} has recently used the derivative process to study convergence rates for RBMs satisfying certain strong uniformity conditions in dimension (which do not hold for the symmetric Atlas model). 

Although currently there are no available lower bounds, it seems  plausible (based on informal calculations using the associated killed Markov chain) that the true $L^1$-Wasserstein distance indeed decays polynomially for $t \le t(d)$, for some $t(d) \rightarrow \infty$ as $d \rightarrow \infty$. The polynomial rates of convergence to stationarity obtained in \cite{banerjee2020rates} for the Potlatch process on $\mathbb{Z}^k$, which (for $k=1$) can be loosely thought of as a `Poissonian version' of the gap process of the infinite Atlas model constructed in \cite{pitman_pal}, lend further evidence to this belief. 

Moreover, dimension-free local convergence is not expected to hold for initial distributions `too far away' from $\ZZ(\infty)$. This belief stems from the analogy with the \emph{infinite} Atlas model (see Section \ref{infatsec}) which has uncountably  many product form stationary distributions $\pi_a := \bigotimes_{i=1}^{\infty} \operatorname{Exp}(2 + ia), \ a \ge 0$ \cite{sarantsev2017stationary2}, only one of which, namely $\pi_0$, is a weak limit of the stationary distribution \eqref{eqn:atlas_stationary} of the $d$-dimensional Atlas model (extended to a measure on $\mathbb{R}^{\infty}_+$ in an obvious manner) as $d \rightarrow \infty$. This leads to the heuristic that, for large $d$, the $d$-dimensional gap process with initial distribution `close' to the projection (onto the first $d$ co-ordinates) of $\pi_a$ for some $a>0$ spends a long time near this projection (which diverges as $d \rightarrow \infty$) before converging to \eqref{eqn:atlas_stationary}. From this heuristic, one expects that dimension-free convergence rates for associated statistics can only be obtained if the initial gap distribution is `close' to the stationary measure \eqref{eqn:atlas_stationary} in a certain sense.

Evidence for this heuristic is provided in the few available results on `uniform in dimension' convergence rates for rank-based diffusions \cite{jourdain2008propagation,jourdain2013propagation}. In both these papers, under strong convexity assumptions on the drifts of the particles, dimension-free exponential ergodicity was proven for the joint density of the particle system when the initial distribution is close to the stationary distribution as quantified by the Dirichlet energy functional (see \cite[Theorem 2.12]{jourdain2008propagation} and \cite[Corollary 3.8]{jourdain2013propagation}). The symmetric Atlas model lacks such convexity in drift and hence, the dimension-free Poincar\'e inequality for the stationary density, that is crucial to the methods of \cite{jourdain2008propagation,jourdain2013propagation}, does not apply.

\section{Domains of Attraction of the infinite Atlas model}\label{infatsec}

In the remaining two sections of this review we will discuss infinite dimensional reflected Brownian motions.
The processes we consider arise from the study of the infinite Atlas model and its various rank based diffusion model variants. This section will consider the standard Atlas model whereas in the next section we consider a more general family of rank based diffusions. 

The Atlas model describes the evolution of a countable collection of particles on the real line moving according to mutually independent Brownian motions and such that, at any point of time, the lowest particle has a unit upward drift.
Letting $\{Y_i(\cdot) : i \in \mathbb{N}_0\}$ denote the state processes of the particles in the system, one can formally describe the evolution by the following system of SDE:
\begin{equation}\label{Atlasdef}
dY_i(s) = \mathbf{1}(Y_i(s) = Y_{(0)}(s))ds + d W_i(s), \ s \ge 0,  i \in \mathbb{N}_0,
\end{equation}
where $Y_{(0)}(t)$ is the state of the lowest particle at time $t$.
Define
\begin{equation}\label{eq:cludef}
	\clu := \Big\{\ybd = (y_0, y_1, \ldots)' \in \RR^{\infty}: \sum_{i=0}^{\infty} e^{-\alpha [(y_i)_+]^2} < \infty \mbox{ for all } \alpha >0\Big\}.
\end{equation}
Following \cite{AS}, a $\ybd \in \RR^{\infty}$ is called {\em rankable} if there is a one-to-one map $\pbd$ from $\NN_0$ onto itself such that $y_{\pbd(i)} \le y_{\pbd(j)}$ whenever $i<j$, $i,j \in \NN_0$. If $\ybd$ is rankable, we will write $y_{(i)} := y_{\pbd(i)}, \, i \in \NN_0$. It is easily seen that any $\ybd \in \clu$ is locally finite and hence rankable. It was shown in \cite[Theorem 3.2 and Lemma 3.4]{AS} that, if the distribution of the initial vector $\YY(0)$ of the particles in the infinite Atlas model is supported within the set $\clu$,
then the system \eqref{Atlasdef} has a weak solution that is unique in law and the solution process $\YY(t) \in \clu$ a.s. for all $t\ge 0$. Consequently, it can be shown that, almost surely, the ranking $(Y_0(t),Y_1(t),\dots) \mapsto (Y_{(0)}(t),Y_{(1)}(t),\dots)$ is well defined (where $Y_{(i)}(t)$ is the state of the $i+1$-th particle from below) for all $t \ge 0$. We will assume throughout that $0 \le Y_{0}(0) \le Y_{1}(0) \le Y_{2}(0) \le \dots$. 
By calculations based on Tanaka's formula (see \cite[Lemma 3.5]{AS}), it follows that the processes defined by
\begin{equation}\label{eq:bistar}
B^*_i(t) := \sum_{j=0}^{\infty}\int_0^t\mathbf{1}(Y_j(s) = Y_{(i)}(s)) dW_j(s), \ i \in \mathbb{N}_0, t \ge 0,
\end{equation}
are independent standard Brownian motions which can be used to write down the following stochastic differential equation for $\{Y_{(i)}(\cdot) : i \in \mathbb{N}_0\}$:
\begin{equation}\label{order_SDE}
dY_{(i)}(t) = \mathbf{1}(i=0) dt + dB^*_i(t) - \frac{1}{2}dL^*_{i+1}(t) + \frac{1}{2} dL^*_{i}(t), \ t \ge 0, i \in \mathbb{N}_0.
\end{equation}
Here, $L^*_{0}(\cdot) \equiv 0$ and for $i \in \mathbb{N}$, $L^*_i(\cdot)$ denotes the local time of collision between the $(i-1)$-th and $i$-th particle.

The gap process in the infinite Atlas model is the $\R_+^{\infty}$-valued process $\Z(\cdot)=$ $(Z_1(\cdot), Z_2(\cdot), \dots)$ defined by
\begin{equation}\label{gapdef}
Z_i(\cdot) := Y_{(i)}(\cdot) - Y_{(i-1)}(\cdot), \ i \in \mathbb{N}.
\end{equation}
We will be primarily interested in the long time behavior of the gap process. 

It was shown in \cite{pitman_pal} that the distribution $\pi := \otimes_{i=1}^{\infty} \operatorname{Exp}(2)$ is a stationary distribution of the gap process $\Z(\cdot)$. It was also conjectured there that this is the unique stationary distribution of the gap process. Surprisingly, this was shown to be not true in \cite{sarantsev2017stationary2} who gave an uncountable collection of product form stationary distributions defined as
\begin{equation}\label{altdef}
\pi_a := \bigotimes_{i=1}^{\infty} \operatorname{Exp}(2 + ia), \ a \ge 0.
\end{equation}
\textcolor{black}{The `maximal' stationary distribution $\pi_0$ (in the sense of stochastic domination) is somewhat special in this collection as described below. For this reason, in what follows, we will occasionally using the notation $\pi$ for the measure $\pi_0$.}

For initial gap distribution $\nu$, we will denote by $\hat \nu_t$ the probability law of $\Z(t)$. Obviously, when $\nu = \pi_a$, $\hat\nu_t  = \pi_a$ for all $t\ge 0$. In general, if $\hat \nu_t$ converges weakly to one of the stationary measures $\pi_a$, $a\ge 0$, we say that $\nu$ is in the {\em Domain of Attraction} (DoA) of $\pi_a$.
Characterizing the domain of attraction of the above collection of stationary measures is a challenging open problem. \textcolor{black}{It was shown in \cite[Theorem 4.7]{AS} using comparison techniques with finite-dimensional  Atlas models that if the law $\nu$ of the  initial gaps $\Z(0)$ stochastically dominates $\pi = \pi_0$ (in a coordinate-wise sense), then $\nu$ is in DoA of $\pi$, i.e. the law of $\Z(t)$ converges weakly to $\pi$ as $t \rightarrow \infty$.} Recently, \cite{DJO} established a significantly larger domain of attraction for $\pi$ using relative entropy and Dirichlet form techniques. They showed that a $\nu$ is in DoA of $\pi$ if the random vector $\Z(0)$ with distribution $\nu$ almost surely satisfies the following conditions for some $\beta \in [1,2)$ and eventually non-decreasing sequence $\{\theta(m) : m \ge 1\}$ with $\inf_m\{\theta(m-1)/\theta(m)\}>0$:
\begin{align}
\label{d1} &\limsup_{m \rightarrow \infty}\frac{1}{m^{\beta}\theta(m)} \sum_{j=1}^m Z_j(0) < \infty,\\
\label{d2} &\limsup_{m \rightarrow \infty}\frac{1}{m^{\beta}\theta(m)} \sum_{j=1}^m (\log Z_j(0))_- < \infty,\\ 
\label{d3} &\liminf_{m \rightarrow \infty}\frac{1}{m^{\beta^2/(1+\beta)}\theta(m)} \sum_{j=1}^m Z_j(0) = \infty,
\end{align}
with the additional requirement that $\theta(m) \ge \log m, m \ge 1,$ if $\beta = 1$. This is a significant advance in our understanding of the domain of attraction properties of the infinite Atlas model. However, the conditions \eqref{d1}-\eqref{d3} involve upper and lower bounding the growth rate of the starting points $Y_m(0) = \sum_{j=1}^m Z_j(0)$ with respect to $m$ in terms of a common parameter $\beta \in [1,2)$, which partially restricts its applicability. In particular, they do not cover all initial gap distributions that stochastically dominate $\pi$ (e.g. $Z_j(0)\sim e^{j^2}$), which were shown to be in DoA of $\pi$ in \cite{AS}. Moreover, one expects that the properties of the first few gaps should not impact the domain of attraction property of a given vector of initial gaps, however \eqref{d2} is not satisfied  if even one of the gaps is zero, which makes this condition somewhat unnatural. This condition arises from the relative entropy methods used in \cite{DJO} (see e.g. the estimate  (3.5) therein) which make an important use of the fact that  the gaps are non-zero.

Beyond these results, little is known about the domain of attraction of stationary measures of the infinite Atlas model. Especially, for the other stationary measures $\pi_a, a >0$, there are no results about the domain of attraction. The fundamental challenge in investigating the latter question is that stability for initial distribution $\pi_a$ for any $a>0$ involves a delicate interplay between the drift of the lowest ranked particle and the `downward pressure' from the cloud of higher particles which has an exponentially increasing density. In particular, one cannot use approximation techniques based on finite dimensional systems as was done in previous works (namely, \cite{AS}  and \cite{DJO}). This can also be observed, as noted previously, upon recalling that the marginals of the stationary distribution of the $d$-dimensional Atlas model in \eqref{eqn:atlas_stationary} approach those of $\pi_0$ as $d \rightarrow \infty$
and so a finite dimensional approximation approach is designed to {\em select} this particular stationary distribution of the infinite Atlas model.

In \cite{banerjee2022domains}, we consider a somewhat weaker formulation of domain of attraction of the stationary measures $\pi_a$, $a \ge 0$.
Specifically, for initial gap distribution $\nu$, define for $t>0$
$$\nu_t := \frac{1}{t} \int_0^t \hat \nu_s ds,$$
where $\hat \nu_t$ is as introduced in the paragraph below \eqref{altdef}.
We say that a $\nu$ is in the {\em Weak  Domain of Attraction} (WDoA) of $\pi_a$ for $a\ge 0$, if $\nu_t \to \pi_a$ weakly as $t\to \infty$. Note that, although the convergence of the time-averaged laws $\nu_t$  to $\pi_a$ (ergodic limit) is implied by the convergence of $\hat \nu_t$ to $\pi_a$ (marginal time limit), the converse is not clear.
However, if $\nu$ is in the WDoA of $\pi_a$ for some $a>0$, it is \emph{necessarily not} in the DoA of $\pi=\pi_0$.  The paper \cite{banerjee2022domains}  provides sufficient conditions for a measure $\nu$ to be in the WDoA of $\pi_a$ for a general $a\ge 0$. The main results of \cite{banerjee2022domains} are described below.

\subsection{Main results}\label{main}
Denote by $\cls_0$ the class of probability measures on $\mathbb{R}_+^{\infty}$ such that the corresponding $\mathbb{R}_+^{\infty}$-valued random variable $\YY(0) = (Y_i(0))_{i\in \NN_0}$ is supported on $\clu$ and 
$$0 = Y_0(0)\le Y_1(0) \le \cdots \mbox{ a.s. }$$
Given a measure $\gamma \in \cls_0$, from weak existence and uniqueness, we can construct a filtered probability space $(\Om, \clf, \Pmb, \{\clf_t\}_{t\ge 0})$ on which are given mutually independent $\clf_t$-Brownian motions $\{W_i(\cdot), i \in \NN_0\}$ and $\clf_t$ adapted continuous processes $\{Y_i(t), t \ge 0, i \in \NN_0\}$ such that $Y_i$ solve \eqref{Atlasdef} with $\Pmb\circ(\YY(0))^{-1} = \gamma$.
On this space the processes $Y_{(i)}$ and $Z_i$ are well defined and the distribution $\Pmb\circ (\Z(\cdot))^{-1}$ is uniquely determined from $\gamma$. Furthermore $\Z$ is a Markov process with values in (a subset of) $\mathbb{R}_+^{\infty}$. Let
$$\cls = \{\Pmb\circ \Z(0)^{-1}: \Pmb \circ \YY(0)^{-1} \in \cls_0\}.$$
Note that $\mu \in \mathcal{S}$ if for some $\mathbb{R}_+^{\infty}$-valued random variable $\YY(0)$ with probability law in $\cls_0$,  the vector $(Y_1(0) - Y_0(0), Y_2(0) - Y_0(0),\dots)$ has  probability law $\mu$.

The first theorem gives a sufficient condition for a probability measure on $\R_+^{\infty}$ to be in the WDoA of $\pi$.
 \begin{theorem}[\cite{banerjee2022domains}]\label{piconv}
 Suppose that the probability measure
 $\mu$ on $\R_+^{\infty}$ satisfies the following: there exists a coupling $(\U,\V)$ of $\mu$ and $\pi$ such that, almost surely,
 \begin{equation}\label{star}
 \liminf_{d \rightarrow \infty} \frac{1}{\sqrt{d} (\log d)} \sum_{i=1}^d U_i \wedge V_i = \infty.
 \end{equation}
 Then $\mu \in \cls$ and it belongs to the WDoA of $\pi$.
 \end{theorem}

 The following corollary gives some natural situations in which \eqref{star} holds.
 
 \begin{corollary}[\cite{banerjee2022domains}]\label{starcheck}
Let $\mu$ be a probability measure on $\R_+^{\infty}$ and
  $\U \sim \mu$. Suppose one of the following conditions hold:
 \begin{itemize}
 \item[(i)] $\mu$ is stochastically dominated by $\pi$ and, almost surely,
 \begin{equation}\label{star1}
 \liminf_{d \rightarrow \infty} \frac{1}{\sqrt{d} (\log d)} \sum_{i=1}^d U_i = \infty.
 \end{equation}
 \item[(ii)] For some a.s. finite random variable $M$
 \begin{equation}\label{star2}
 \liminf_{d \rightarrow \infty} \frac{1}{\sqrt{d} (\log d)} \sum_{i=1}^d (U_i \wedge M) = \infty, \mbox{ a.s. }
 \end{equation}
 \item[(iii)] $U_i = \lambda_i \Theta_i, \ i \in \mathbb{N}$, where $\{\Theta_i\}$ are iid non-negative random variables satisfying $\mathbb{P}(\Theta_1>0)>0$, and \textcolor{black}{$\{\lambda_i\}$ are positive deterministic real numbers satisfying one of the following: 
 \begin{itemize}
 \item[(a)] $\liminf_{j \rightarrow \infty} \lambda_j>0$,\\
 \item[(b)] $\limsup_{j \rightarrow \infty}\lambda_j < \infty$ and
 \begin{equation}\label{star3}
 \liminf_{d \rightarrow \infty}\frac{1}{\sqrt{d} (\log d)} \sum_{i=1}^d \lambda_i = \infty.
 \end{equation}
 \end{itemize}}
 \end{itemize}
 Then \eqref{star} holds. Hence, in all the above cases, $\mu \in \cls$ and is in the WDoA of $\pi$.

 \end{corollary}

 \begin{remark}
	 \label{remthm1}
	 We make the following observations.
	 \begin{enumerate}[(a)]
	 	\item Suppose that $\pi$ is stochastically dominated by $\mu$. Then \eqref{star} holds by the strong law of large numbers, and so Theorem \ref{piconv} shows that $\mu \in \cls$ and is in the WDoA of $\pi$. \textcolor{black}{In fact in this case  \cite[Theorem 4.7]{AS} shows that $\mu$ is in the DoA of $\pi$.}
		\item Note that for any $C \in (0,\infty)$, 
		$$((2C) \wedge 1) [U_i \wedge (1/2)] \le U_i \wedge C \le ((2C) \vee 1) [U_i \wedge (1/2)].$$
		 Hence, if \eqref{star2} holds for one a.s. finite random variable $M$ then it holds for every such random variable.
		 In particular if $\limsup_{n\to \infty} U_n <\infty$ a.s. and \eqref{star1} is satisfied then \eqref{star} holds and so $\mu \in \cls$ and is in the WDoA of $\pi$.
		 \item The paper \cite{DJO} notes two important settings where conditions \eqref{d1}-\eqref{d3} are satisfied. These are: (i) for some $c \in [1,\infty)$, $\lambda_j \in [c^{-1}, c]$ for all $j\in \NN$ and 
		 $ Z_j(0) = \lambda_j \Theta_j$ where $\Theta_j$ are iid with finite mean and such that 
		 $E\log (\Theta_1)_ {-}<\infty$, and (ii) $ Z_j(0) = \lambda_j \Theta_j$ where $\Theta_j$ are iid exponential  with mean $1$ and, either $\lambda_d\downarrow 0$ and $\frac{1}{\sqrt{d}\log d} \sum_{i=1}^d \lambda_i \to \infty$, as $d\to \infty$, or $\lambda_d\uparrow \infty$ and $\limsup_{d\to \infty}\frac{1}{d^{\beta}} \sum_{i=1}^d \lambda_i < \infty$ for some $\beta <2$. We note that conditions assumed in the above settings are substantially stronger than the one assumed in part (iii) of Corollary \ref{starcheck}. \textcolor{black}{However  \cite{DJO}  establishes the stronger  marginal time convergence as opposed to a time averaged limit considered in Corollary \ref{starcheck}. See part (g) on a related open problem.}
		 \item In the setting of (iii) the sequence $\{\lambda_i\}$  satisfies \eqref{star3} if 
		  $(\log i)/\sqrt{i} \lambda_i \to 0$ as $i \to \infty$.
		 Indeed, note that, for any $C>0$, there exists $i_C \in \mathbb{N}$ such that for all $i \ge i_C$,
		  $$
		  \sum_{j=i_C}^i \lambda_j \ge C\sum_{j=i_C}^i\frac{\log j}{\sqrt{j}} \ge C\int_{i_C}^{i+1}\frac{\log z}{\sqrt{z}}dz.
		  $$
		  \textcolor{black}{Using the fact that the primitive of $\log z/\sqrt{z}$ is $2\sqrt{z}(\log z - 2)$, $\{\lambda_j\}$ is seen to satisfy \eqref{star3}}
		  and so, if in addition $\limsup_{j \rightarrow \infty}\lambda_j < \infty$, we have that
		  the conditions of (iii) (b) are satisfied.
		 \item  We note that the conditions prescribed in Theorem \ref{piconv} and Corollary \ref{starcheck} all involve quantifying how close together the particles can be on the average in the initial configuration of the Atlas model. In particular, we do not require any 
		 upper bounds on the sizes of $U_i$.
		 This is in contrast with condition \eqref{d1} assumed in \cite{DJO} which requires the particles to be not too far apart on the average.  Intuitively one expects convergence to $\pi$ to hold when initially the particles are not too densely packed and upper bounds on the  average rate of growth of the initial spacings are somewhat unnatural. The result in \cite[Theorem 4.7]{AS} also points to this heuristic by showing weak convergence to $\pi$ from all initial gap configurations that stochastically dominate $\pi$. 
		 We also note that Theorem \ref{piconv} and Corollary \ref{starcheck}
		  allow the first few gaps of the initial gap sequence to be arbitrary and do not require any condition analogous to \eqref{d2} and can, in particular, allow gaps to be zero.

		\item  We note that the conditions \eqref{d1}-\eqref{d3} of \cite{DJO} do not  imply \eqref{star}. To see this, consider the deterministic sequence of initial gaps: $U_i = i^{-2/3}$ for $n^3 < i < (n+1)^3$ and $U_{n^3} = n$, for any $n \in \mathbb{N}$. It can be checked that, $\sum_{i=1}^d U_i$ grows like $d^{2/3}$ and $\sum_{i=1}^d (\log U_i)_-$ grows like $d\log d$ as $d \rightarrow \infty$. Thus,  \eqref{d1}-\eqref{d3} of \cite{DJO} hold with $\beta=1$ and $\theta(d)= \log d$. However, almost surely, $\sum_{i=1}^d U_i \wedge V_i$ grows like $d^{1/3}$, and therefore, \eqref{star} is violated.
		\item It will be interesting to investigate whether the condition in Theorem \ref{piconv} in fact  implies the stronger property  that $\mu$ is in the DoA of $\pi$. We leave this as an open problem.
	 \end{enumerate}
 
 \end{remark}
 
The next theorem provides sufficient conditions for a measure on $\R_+^{\infty}$ to be in the WDoA of
 one of the other stationary measures $\pi_a, a>0$.
 
 \begin{theorem}[\cite{banerjee2022domains}]\label{piaconv}
 Fix $a>0$. Suppose that the probability measure $\mu$ on $\R_+^{\infty}$ satisfies the following:
 There exists a coupling $(\U, \V_a)$ of  $\mu$  and $\pi_a$ such that, almost surely,
\begin{equation}\label{stara}
 \limsup_{d \rightarrow \infty} \frac{\log \log d}{\log d} \sum_{i=1}^d |V_{a,i} - U_i| = 0, \mbox{ and } \limsup_{d\to \infty} \frac{U_d}{dV_{a,d}}<\infty.
 \end{equation}
 Then $\mu \in \cls$ and it belongs to the WDoA of $\pi_a$.
 \end{theorem}
\textcolor{black}{
We remark that Theorem \ref{piaconv} also holds for the case $a=0$. Indeed, suppose that
the condition \eqref{stara} in the theorem holds for $a=0$. Denoting the corresponding  $\V_0$ as $\V$, we have that
 $$
 \sum_{i=1}^d U_i \wedge V_i = \sum_{i=1}^d U_i \vee V_i - \sum_{i=1}^d |V_i - U_i| \ge \sum_{i=1}^d V_i - \sum_{i=1}^d |V_i - U_i|.
 $$
 Using this and the law of large numbers for $\{V_i\}$, it follows that 
  \eqref{star} holds and thus from Theorem \ref{piconv} the conclusion of Theorem \ref{piaconv} holds with $a=0$. The reason we do not note the case $a=0$ in the statement of Theorem \ref{piaconv} is because \eqref{star} is a strictly weaker assumption than \eqref{stara} when $a=0$. This is expected as for any initial gap distribution that stochastically dominates $\pi$, convergence to $\pi$ holds (see Remark \ref{remthm1} (a)). However, for convergence to $\pi_a$ for some $a>0$, the initial gap distribution should be appropriately close to $\pi_a$ (see Remark \ref{remthm2} (b) below).}
 
 The following corollary gives a set of random initial conditions for which \eqref{stara} holds.
 
 \begin{corollary}\label{staracheck}
 Fix $a>0$. Let
 \begin{equation}\label{eq:805nn}
 \mu := \bigotimes_{i=1}^{\infty} \operatorname{Exp}(2 + ia + \lambda_i),
 \end{equation}
 where $\{\lambda_i\}$ are real numbers satisfying $\lambda_i > -(2+ia)$ for all $i$ and
 \begin{equation}\label{eq:802n}
 	-a < \liminf_{i\to \infty} \frac{\lambda_i}{i} \le \limsup_{i\to \infty} \frac{\lambda_i}{i} < \infty.
 \end{equation}
 %
%
Moreover, assume
 \begin{equation}\label{aexp}
 \limsup_{d \rightarrow \infty} \frac{\log \log d}{\log d}\sum_{i=1}^d\frac{|\lambda_i|}{i^2} = 0. 
 \end{equation}
 Then \eqref{stara} holds, and hence, $\mu \in \cls$ and belongs to the WDoA of $\pi_a$. 
 \end{corollary}
 
 \begin{remark}
	 \label{remthm2}
	 We note the following.
	 \begin{enumerate}[(a)]
		 \item Clearly, there are measures $\mu$ of the form in \eqref{eq:805nn}
		  for which $\limsup_{i\to \infty} \frac{\lambda_i}{i} < \infty$,
		  $\liminf_{i\to \infty}\frac{|\lambda_i|}{i}>0$ (and so \eqref{eq:802n} holds) and which do not lie in the WDoA of $\pi_a$ (for example, $\pi_{a'}$ for any $a' \neq a$). However, the corollary says that if $|\lambda_i|$ grows slightly slower than $i$ then $\mu$ is indeed in the WDoA of $\pi_a$. More precisely, the convergence  in \eqref{aexp} holds in particular if 
		  $\limsup_{i\to \infty} |\lambda_i|\log\log i/i=0$.
		  %
		  %
	Indeed,  note that for any $\delta>0$, there exists sufficiently large $i_\delta \in \mathbb{N}$ such that for $i \ge i_\delta$,
\begin{equation}\label{ac2}
\sum_{j=i_\delta}^i \frac{|\lambda_j|}{j^2} \le \sum_{j=i_\delta}^i \frac{\delta_j}{j \log \log j} \le \delta\int_{i_\delta-1}^{i} \frac{1}{z\log \log z}dz = \delta\int_{\log(i_\delta-1)}^{\log i} \frac{1}{\log w}dw.
\end{equation}
As $w \mapsto \frac{1}{\log w}$ is a slowly varying function, by Karamata's theorem \cite[Theorem 1.2.6 (a)]{mikosch1999regular}, 
$$
 \lim_{i \rightarrow \infty} \frac{\log \log i}{\log i}\int_{\log(i_0-1)}^{\log i} \frac{1}{\log w}dw = 1.
 $$
 Thus, we conclude from \eqref{ac2} that for any $\delta>0$,
 $$
  \limsup_{i \rightarrow \infty} \frac{\log \log i}{\log i} \sum_{j=1}^i \frac{|\lambda_j|}{j^2} \le \delta.
 $$
As $\delta>0$ is arbitrary, \eqref{aexp} holds.

 \item Note that if $\U \sim \pi$, then for $a, a'>0$,
 $$V_{a,i} = \frac{2}{2+ia} U_i, \; V_{a',i} = \frac{2}{2+ia'} U_i, \; i \in \NN$$
 defines a coupling of $\pi_a$ and $\pi_{a'}$.
 For this coupling,
 \begin{equation}\label{eq:203}
	 \sum_{i=1}^d |V_{a,i}- V_{a',i}| \sim O(|a-a'| \log d) \mbox{ and } \limsup_{d\to \infty} \frac{V_{a', d}}{dV_{a,d}} = 0.
	 \end{equation}
Since $\pi_{a'}$ is obviously not in WDoA of $\pi_a$ for $a\neq a'$, the first
 property in  \eqref{eq:203} says that the first requirement in \eqref{stara} is not far from what one expects.
  %
  %
  We conjecture that for $\mu$ to be in the WDoA of $\pi_a$  it is necessary that $(\log d)^{-1}\sum_{i=1}^d |V_{a,i} - U_i| \rightarrow 0$ as $d \rightarrow \infty$ 
 for some coupling $(\U, \V_a)$ of  $\mu$  and $\pi_a$.  
 Theorem \ref{piaconv}  says that if the above convergence to $0$ holds at a rate faster than $1/(\log \log d)$ for some coupling of $\mu$  and $\pi_a$ then that (together with the second condition in \eqref{stara}) is sufficient for 
 $\mu$ to be in the WDoA of $\pi_a$.
 
Further, observe that if the second condition in \eqref{stara} does not hold, then there are infinitely many $d \in \mathbb{N}$ such that the lowest $d+1$ particles are separated from the rest by a large  gap. If one such gap is very large, then it could happen that the distribution of gaps between the $d+1$ particles stabilizes towards the unique stationary gap distribution $\pi^{(d)}$ of the corresponding finite Atlas model before the remaining particles have interacted with them. 
As $\pi^{(d)}$ converges weakly to $\pi$ as $d \rightarrow \infty$, one does not expect in such situations the convergence (of time averaged laws) to $\pi_a$
to hold for any $a>0$.
  \end{enumerate}
 \end{remark}

\subsection{Outline of Approach}
Proofs of the above results are based on a pathwise approach to studying the long time behavior of the infinite Atlas model using synchronous couplings, namely, two versions of the infinite Atlas model described in terms of the ordered particles via \eqref{order_SDE} and started from different initial conditions but driven by the same collection of Brownian motions. Central ingredients in the proofs are suitable estimates on the decay rate of the $L^1$ distance between the gaps in synchronously coupled  ordered infinite Atlas models. These estimates are obtained by analyzing certain excursions 
of the difference of the coupled processes where each excursion ensures the contraction of the $L^1$ distance by a fixed deterministic amount. The key is to appropriately control the number and the lengths of such excursions. Certain monotonicity properties of synchronous couplings and a quantification of the influence of far away coordinates on the first few gaps also play a crucial role. Using these tools, the discrepancy between the gap processes of two synchronously coupled infinite Atlas models can be controlled in terms of associated gap processes when the starting configurations differ only in finitely many coordinates. 
The latter are more convenient to work with as for them the 
initial $L^1$ distance between the  gap processes is finite, and the aforementioned excursion analysis can be applied to obtain the main results.

\section{Extremal and product form 
invariant distributions of infinite rank-based diffusions}
\label{sec:extrem}

The results described in the last section give us some understanding about the local stability structure of the set of invariant distributions of the infinite Atlas model. However, it does not say anything about the geometry of the (convex) set of invariant measures. In particular, a longstanding open problem is to characterize all extremal invariant measures (corner points of the convex set). In this section, we outline some recent results from \cite{banerjee2022extremal}, which makes progress in this direction.

We will in fact consider the more general setting  of an infinite rank-based diffusion where the particles move on the real line according to mutually independent Brownian motions, with the $i$-th particle from the left getting a constant drift of $g_{i-1}$. The vector $\gbd := (g_0,g_1,\dots)'$ is called the drift vector and we call this model the \emph{$\gbd$-Atlas model}.
Such particle systems were originally introduced in stochastic portfolio theory \cite{Fern,fernholz2009stochastic} as models for stock growth evolution in equity markets and have been investigated extensively in recent years in several different directions. In particular, characterizations of such particle systems as scaling limits of jump processes with local interactions on integer lattices, such as the totally asymmetric simple exclusion process, have been studied in \cite{karatzas_skewatlas}. Various types of results for  the asymptotic behavior of the empirical measure of the particle states have been studied, such as propagation of chaos, characterization of the associated McKean-Vlasov equation and nonlinear Markov processes 
\cite{shkol, jourdain2008propagation}, large deviation principles \cite{dembo2016large}, characterizing the asymptotic density profile and the  trajectory of the leftmost particle via Stefan free-boundary problems \cite{cabezas2019brownian} .  These particle systems also have  close connections with Aldous' ``Up the river'' stochastic control problem \cite{aldous41up}, recently solved in \cite{tang2018optimal}.
Results on wellposedness of the associated stochastic differential equations (in the weak and strong sense) and on absence of triple collisions (three particles at the same place at the same time) have been studied in \cite{bass1987uniqueness, ichiba2013strong, ichiba2010collisions,sarantsev2015triple,ichiba2017yet}.

As in the case of the infinite Atlas model, the natural state space of the location of the particles in the more general system is also given by $\clu$ defined in \eqref{eq:cludef}.
For a sequence $\{W_i\}_{i\in \NN_0}$ of mutually independent standard Brownian motions given on some filtered probability space, and $\ybd = (y_0, y_1, \ldots)' \in \clu$, $\gbd = (g_0, g_1, \ldots)' \in \cld$, consider the following system of equations.
\begin{equation}\label{eq:unrank}
	dY_i(t) = \left[\sum_{k=0}^{\infty} \one(Y_i(t) = Y_{(k)}(t)) g_k\right] dt + dW_i(t), \ t \ge 0,
\end{equation}
with $Y_i(0) = y_i$,  $i \in \NN_0$. Write $\Ybd(\cdot) := (Y_0(\cdot),Y_1(\cdot),\dots)'$.

It has been shown in \cite{AS} that if $\gbd \in \cld$, where
\begin{equation}\label{eq:clddefn}
	\cld := \Big\{\gbd = (g_0, g_1, \ldots)' \in \RR^{\infty}:  \sum_{i=0}^{\infty} g_i^2 <\infty \Big\},
\end{equation}
then for any $\ybd \in \clu$, $P(\Ybf(t)\in \clu \mbox{ for all } t \ge 0)=1$, and there is a unique weak solution $\Ybd(\cdot)$ to \eqref{eq:unrank}.
Analogous to the standard Atlas model setting, the ordered process $\{Y_{(i)}(\cdot) : i \in \NN_0\}$ can be described by the SDE:
\begin{equation}\label{order_SDEg}
dY_{(i)}(t) = g_i dt + dB^*_i(t) - \frac{1}{2}dL^*_{i+1}(t) + \frac{1}{2} dL^*_{i}(t), \ t \ge 0, \;\;   Y_{(i)}(0) = y_{(i)},\, i \in \mathbb{N}_0.
\end{equation}
where $\{B_i^*(\cdot) : i \in \NN_0\}$ are independent Brownian motions defined in \eqref{eq:bistar}, $L^*_{0}(\cdot) \equiv 0$ and for $i \in \mathbb{N}$, $L^*_i(\cdot)$ denotes the local time of collision between the $(i-1)$-th and $i$-th particles.
The gap process for the  $(\gbd, \ybd)$-infinite Atlas model is the $\RR_+^{\infty}$-valued process $\Zbd(\cdot)=$ $(Z_1(\cdot), Z_2(\cdot), \dots)'$ defined by
\begin{equation}\label{gapdefg}
Z_i(\cdot) := Y_{(i)}(\cdot) - Y_{(i-1)}(\cdot), \ i \in \mathbb{N}.
\end{equation}
The natural state space for the gap process is given by
\begin{equation}\label{clvdef}
\clv := \{\zbd \in \RR_+^{\infty}: \mbox{ for some } \ybd \in \clu, \zbd = (y_{(1)}-y_{(0)}, y_{(2)}-y_{(1)}, \ldots)'\}.
\end{equation}

For this setting it is known from the work of \cite{sarantsev2017stationary2} that, once more, the process  associated with the gap sequence of the ranked particle system has a continuum of product-form stationary distributions given as
\begin{equation}\label{eq:piagbddefn}
	\pi_a^{\gbd} := \otimes_{n=1}^{\infty} \mbox{Exp}(n (2\bar g_n+a)), \;\; a > -2 \inf_{n \in \NN} \bar g_n,\end{equation}
where $\bar g_n := \frac{1}{n}(g_0 + \cdots + g_{n-1})$.  In the special case where $\gbd \in \cld_1$, where
\begin{equation}\label{cldzdefn}
	\cld_1 := \{\gbd \in \cld: \text{ there exist } N_1 < N_2 < \dots \rightarrow \infty \text{ such that } \bar g_k > \bar g_{N_j}, \ k=1,\dots,N_j-1, \text{ for all } j \ge 1\},\end{equation}
 $\pi_a^{\gbd}$ is also an invariant distribution for $a=-2 \inf_{n \in \NN} \bar g_n  = -2\lim_{j \rightarrow \infty} \bar g_{N_j}$; in particular for the infinite Atlas model $\pi_a = \pi_a^{\gbd^1}$ is invariant for all $a\ge 0$.

Using Kakutani's theorem, it is easy to verify  that, for different values of $a$, the probability measures 
$\pi_a^{\gbd}$ are mutually singular. These distributions are also special in that they have a product form structure. In particular, if the initial distribution of the gap process is chosen to be one of these distributions, then the laws of distinct gaps at any fixed time are independent despite these gaps having a highly correlated temporal evolution mechanism (see \eqref{order_SDEg}-\eqref{gapdefg}).

Now we present our main results from \cite{banerjee2022extremal} which show that the above distributions are also extremal. Moreover, in some sense, these are the only product form invariant distributions.

\subsection{Main results}

We begin with some notation.
Let $\XX :=  \clc([0,\infty): \RR_+^{\infty})$ which is equipped with the topology of local uniform convergence (with $\RR_+^{\infty}$ equipped with the product topology). 
Define measurable maps $\{\theta_t\}_{t\ge 0}$ from  $(\XX, \mathcal{B}(\XX))$ to itself as
$$\theta_t(\Zbd)(s) := \Zbd(t+s), \; t\ge 0, s\ge 0, \; \Zbd \in \XX.$$
Given $\gbd \in \cld$ and $\zbd \in \clv$, we denote the probability distribution of the gap process of the $\gbd$-Atlas model on $(\XX, \mathcal{B}(\XX))$, with initial gap sequence $\zbd$, by
$\mathbb{P}^{\gbd}_{\zbd}$. Also, for $\gamma \in \clp(\RR_+^{\infty})$ supported on $\clv$ (namely, $\gamma(\clv)=1$), let $\mathbb{P}^{\gbd}_{\gamma} := \int_{\RR_+^{\infty}} \mathbb{P}^{\gbd}_{\zbd} \;\gamma(d\zbd)$. The corresponding expectation operators will be denoted as $\mathbb{E}^{\gbd}_{\zbd}$ and $\mathbb{E}^{\gbd}_{\gamma}$ respectively.
Denote by $\mathcal{I}^{\gbd}$ the collection of all invariant (stationary) probability measures of the gap process of the $\gbd$-Atlas model supported on $\clv$, namely
$$ \mathcal{I}^{\gbd} := \{\gamma \in \clp(\XX):  \gamma(\clv^c)=0, \mbox{ and } \mathbb{P}^{\gbd}_{\gamma} \circ \theta_t^{-1} = \mathbb{P}^{\gbd}_{\gamma} \mbox{ for all } t \ge 0\}.$$
%

Abusing notation, the canonical coordinate process on $(\XX, \mathcal{B}(\XX))$ will be denoted by $\{\Zbd(t)\}_{t\ge 0}$.
 Let $\mathcal{M}(\RR_+^{\infty})$ be the collection of all real measurable maps on $\RR_+^{\infty}$. For $f \in \mathcal{M}(\RR_+^{\infty})$, $\zbd \in \RR_+^{\infty}$ and $t\ge 0$ such that $\mathbb{E}^{\gbd}_{\zbd}(|f(\Zbd(t))|) <\infty$
 we write
 $$T_t^{\gbd}f(\zbd) := \mathbb{E}^{\gbd}_{\zbd}(f(\Zbd(t))).$$
 For $\gamma \in \clp(\RR_+^{\infty} )$
 let $L^2(\gamma)$ be the collection of  all measurable $\psi: \RR_+^{\infty} \to \RR$ such that
$\int_{\RR_+^{\infty}} |\psi(\zbd)|^2  \gamma(d\zbd) <\infty$.
We denote the inner-product and the norm on $L^2(\gamma)$ as $\langle \cdot, \cdot \rangle$ and $\|\cdot\|$ respectively.
Note that for $\gbd \in \cld$, $\gamma \in \mathcal{I}^{\gbd} $, and
$\psi \in L^2(\gamma)$, $T_t^{\gbd}\psi$ is $\gamma$ a.e. well defined and belongs to $L^2(\gamma)$. Furthermore, the collection $\{T_t^{\gbd}\}_{t\ge 0}$ defines a contraction semigroup on $L^2(\gamma)$, namely
$$T_t^{\gbd}T_s^{\gbd} \psi = T_{t+s}^{\gbd}\psi, \mbox{ and } \|T_t^{\gbd}\psi\| \le \|\psi\| \mbox{ for all } s,t \ge 0 \mbox{ and } \psi \in L^2(\gamma).$$

We now recall the definition of extremality and ergodicity. 
Let, for $\gbd$  as above and $\gamma \in \mathcal{I}^{\gbd}$, $\mathbb{I}_{\gamma}^{\gbd} $ be the collection of all $T_t^{\gbd}$-invariant square integrable functions, namely,
$$\mathbb{I}_{\gamma}^{\gbd} := \{\psi \in L^2(\gamma): T_t^{\gbd}\psi = \psi,\;  \gamma  \mbox{ a.s., for all } t \ge 0\}.$$
We denote the projection of a $\psi \in L^2(\gamma)$ on to the closed subspace $\mathbb{I}_{\gamma}^{\gbd}$ as $\hat \psi_{\gamma}^{\gbd}$.
\begin{definition}
	Let $\gbd\in \cld$. A  $\nu \in \mathcal{I}^{\gbd}$ is said to be an extremal invariant distribution of the gap process of the $\gbd$-Atlas model if, whenever for some $\varepsilon \in (0,1)$ and  $\nu_1, \nu_2 \in \mathcal{I}^{\gbd}$  we have $\nu = \varepsilon \nu_1 + (1-\varepsilon)\nu_2$, then $\nu_1 = \nu_2 = \nu$.
	We denote the collection of all such measures by $\mathcal{I}^{\gbd}_e$.
	
	We call $\nu \in \mathcal{I}^{\gbd}$ an ergodic probability measure for the invariant distribution of the gap process of the $\gbd$-Atlas model if for all
	$\psi \in L^2(\nu)$,  $\hat \psi_{\nu}^{\gbd}$ is constant $\nu$-a.s. We denote the collection of all such measures by $\mathcal{I}^{\gbd}_{er}$.
\end{definition}
We note that (cf. proof of Lemma \ref{lemerg} below) if $\gamma \in \mathcal{I}^{\gbd}_{er}$, then for any $\psi \in L^2(\gamma)$, and $\gamma$ a.e. $\zbd$,
$$\frac{1}{t} \int_0^t T_s^{\gbd}\psi(\zbd) ds \to \int_{\RR_+^{\infty}} \psi(\xbd) \gamma(\xbd), \mbox{ in } L^2(\gamma), \mbox{ as } t \to \infty.$$

The following result, which says that extremal invariant measures and ergodic invariant measures are the same, is standard.
\begin{lemma}\label{lemerg}
	Let $\gbd \in \cld$. Then $\mathcal{I}^{\gbd}_e = \mathcal{I}^{\gbd}_{er}$.
	Let $\gamma \in \mathcal{I}^{\gbd}$ and suppose that for every bounded measurable $\psi:\RR_+^{\infty} \to \RR$, 
	$\hat \psi_{\gamma}^{\gbd}$ is constant, $\gamma$ a.s. Then
	$\gamma \in \mathcal{I}^{\gbd}_e$.
\end{lemma}

The following is the first main result of \cite{banerjee2022extremal}.
\begin{theorem}[\cite{banerjee2022extremal}]\label{thm_main}
	Let $\gbd\in \cld$. Then, for every $a > -2 \inf_{n \in \NN} \bar g_n$, $\pi_a^{\gbd} \in \mathcal{I}^{\gbd}_e = \mathcal{I}^{\gbd}_{er}$.  Also, when $\gbd \in \cld_1$, $\pi_a^{\gbd} \in \mathcal{I}^{\gbd}_e = \mathcal{I}^{\gbd}_{er}$ also for $a = -2 \inf_{n \in \NN} \bar g_n$.
\end{theorem}
The above theorem proves the extremality of the invariant measures $\pi_a^{\gbd}$ for suitable values of $a$. One can ask whether these are the only extremal invariant measures of the gap process of the $\gbd$-Atlas model supported on $\clv$. The answer to this question when $\gbd = \mathbf{0}$ is affirmative from results of \cite[Theorem 4.2]{RuzAiz} , if one restricts to extremal measures satisfying certain  integrability constraints on the denseness of particle configurations.
For a general $\gbd \in \cld$ this is currently a challenging open problem. However, the next result  makes partial progress towards this goal  by showing that for any $\gbd \in \cld_1$ (and under a mild integrability condition), the collection $\{\pi_a^{\gbd}\}$ exhausts all the extremal product form invariant distributions. In fact we prove the substantially stronger statement that  the measures $\pi_a^{\gbd}$ are the only product form (extremal or not) invariant distributions under a mild integrability condition. Qualitatively, this result says  that these measures are  the only invariant distributions that preserve independence of the marginal laws of the gaps in time.
\begin{theorem}[\cite{banerjee2022extremal}]\label{thm_main2}
	Let $\gbd \in \cld_1$ and let $\pi \in \mathcal{I}^{\gbd}$ be a product measure, i.e. $\pi = \otimes_{i=1}^{\infty} \pi_i$ for some $\pi_i \in \clp(\RR_+)$, $i \in \NN$. Suppose that
    \begin{equation}\label{eq:satinteg}
        \int_{\RR_+^{\infty}} \left(\sum_{j=1}^{\infty}
        e^{-\frac{1}{4} (\sum_{l=1}^j z_l)^2} \right) \pi(d\zbd) <
        \infty.
    \end{equation}
    Then, for some $a\ge -2\lim_{j \rightarrow \infty} \bar g_{N_j}$, $\pi = \pi_a^{\gbd}$.
\end{theorem}
 Recall that $\clv$ defined in \eqref{clvdef} consists of $\zbd \in \RR_+^{\infty}$ for which $\sum_{j=1}^{\infty}e^{-\alpha (\sum_{l=1}^j z_l)^2} < \infty$ for all $\alpha>0$. In comparison, condition \eqref{eq:satinteg} requires a finite expectation of $\sum_{j=1}^{\infty}e^{-\frac{1}{4}(\sum_{l=1}^j z_l)^2}$ when $\zbd$ is distributed as $\pi$. 
Roughly speaking, condition  \eqref{eq:satinteg}
 puts a restriction on the rate of increase of the density of particles as one moves away from the lowest ranked particle.

\subsection{Connection with prior work in interacting particle systems}
 Questions about extremality and ergodicity of stationary distributions have been addressed previously in the context of interacting particle systems on countably infinite graphs (see \cite{liggett1976coupling,andjel1982invariant,sethuraman2001extremal,balazs2007existence} and references therein). However, in all these cases, the interactions are Poissonian which enables one to use the (explicit) generator of associated processes in an effective way. The interactions in rank based diffusions are very `singular' owing to the local time based dynamics (see \eqref{order_SDE}) and generator based methods seem to be less tractable. Furthermore, unlike previous works,  the state space  for the gap process (i.e. $\mathbb{R}_+^{\infty}$) is not countable and has a non-smooth boundary, and the process has intricate interactions (oblique reflections) with the boundary. Hence, proving extremality requires new techniques. Our proofs are based on constructing appropriate couplings for these infinite dimensional diffusions which then allow us to prove suitable invariance properties and a certain `directional strong Feller property'. Although our setting and methods are very different, at a high level, the approach we take is inspired by the papers \cite{sethuraman2001extremal,balazs2007existence}.
 
Complete characterization of extremal invariant measures is, in general, a very hard problem. This problem has been completely solved in a few examples of interacting particle systems such as the simple exclusion process \cite{liggett1976coupling} and the zero range process \cite{andjel1982invariant} where the extremal probability measures are fully characterized as an explicit collection of certain product form measures. However, in these models the particle density associated with distinct extremal measures are scalar multiples of each other owing to certain homogeneity properties in the dynamics (see, for example, \cite[Theorem 1.10]{andjel1982invariant}). This, along with the Poissonian nature of the interactions, enables one to prove useful monotonicity properties of the `synchronously coupled' dynamics (see \cite[Section 2]{liggett1976coupling} and \cite[Section 4]{andjel1982invariant}) using generator methods that are crucial to the above characterization results. A key challenge in extending these methods to rank based diffusions of the form considered here is that, in addition to the singular local time interactions, the point process associated with the configuration of particles with gaps distributed as $\pi_a^{\gbd}$ has an intensity function  $\rho(x)$, that grows exponentially as $x\to \infty$ when $a>0$ and, due to a nonlinear dependence on $a$, lacks the scalar multiple property for distinct values of $a$. This is a direct consequence of the inhomogeneity of the topological interactions in our particle systems where the local stability in a certain region of the particle cloud is affected both by the density of particles in the neighborhood and their relative ranks in the full system.  
Moreover, unlike the above interacting particle systems,  in rank based diffusions, when the initial gaps are  given by a stationary distribution,  the point process of particles is typically not stationary. This phenomenon, where  the gaps are stationary while the associated point process  is not, referred to as \emph{quasi-stationarity} in \cite{RuzAiz}, is technically challenging. We note that this latter paper studies one setting where the intensity function grows exponentially and the particle density lacks the scalar multiple property for distinct values of $a$. However their setting, in the context of rank based diffusions, corresponds to the
case  $\gbd = (0,0, \ldots)'$, where the unordered particles behave like independent standard Brownian motions, and this fact is crucially exploited in \cite{RuzAiz}.

\subsection{Proof outline}

We now make some comments on proof ideas. The key step in proving the extremality of $\pi_a^{\gbd}$ is to establish that any bounded measurable function $\psi$ on $\RR_+^{\infty}$ that is $\pi_a^{\gbd}$-a.e. invariant, under the action of the semigroup of the Markov process corresponding to the $\gbd$-Atlas gap sequence, is constant $\pi_a^{\gbd}$-a.e.
If $\gbd = \gbd^1$ and $a=0$, we have that $\pi_a^{\gbd} = \otimes_{i=1}^{\infty} \mbox{Exp}(2)$, and therefore the coordinate sequence $\{Z_i\}_{i=1}^{\infty}$ is iid under $\pi_a^{\gbd}$. In this case, from the Hewitt-Savage zero-one law it suffices to show that $\psi$ is $\pi_a^{\gbd}$-a.e. invariant under all finite permutations of the coordinates of $\RR_+^{\infty}$. For this, in turn, it suffices to simply prove the above invariance property for transpositions of the $i$-th and $(i+1)$-th coordinates, for all $i \in \NN$. For a general $\pi_a^{\gbd}$, the situation is a bit more involved as the coordinate sequence $\{Z_i\}_{i=1}^{\infty}$  is not iid any more. Nevertheless, from the scaling properties of Exponential distributions it follows that, with $c_n := 2[n (2\bar g_n+a)]^{-1}$, the sequence 
$\{\tilde Z_n\}_{n\ge 1}$, defined as $\tilde Z_n = c_n^{-1} Z_n$, $n \in \NN$, is iid under $\pi_a^{\gbd}$. In this case, in order to invoke the Hewitt-Savage zero-one law, one needs to argue that for each $i$ the map $\psi$ is $\pi_a^{\gbd}$-a.e. invariant under the transformation that takes the $(i, i+1)$ coordinates $(z_i, z_{i+1})$ to $(\frac{c_i}{c_{i+1}} z_{i+1}, \frac{c_{i+1}}{c_{i}} z_{i})$ and keeps the remaining coordinates the same. Establishing this property is at the heart of the proof of Theorem \ref{thm_main}. A key technical idea in the proof is the construction of a mirror coupling of the first $i+1$ Brownian motions, and synchronous coupling of the remaining Brownian motions,
in the evolution of the ranked infinite $\gbd$-Atlas model corresponding to a pair of nearby initial configurations. Estimates on the probability of a successful coupling, before any of the first $i$-gap processes have hit zero or the lowest $i$-particles have interacted with the higher ranked particles (in a suitable sense), are some of the important ingredients in the proof.

The proof of Theorem \ref{thm_main2} hinges on establishing a key identity for expectations, under the given product form invariant measure $\pi$, of certain integrals involving the state process of the $i$-th gap and the collision local time for the $(j-1)$-th and $j$-th particle, for $i\neq j$. This identity is a consequence of the product form structure of $\pi$ and basic results on local times of continuous semimartingales. 
By using the form of the dynamics of the $\gbd$-Atlas model, the above identity allows us to compute explicitly the Laplace transform of the $i$-th coordinate projection of $\pi$, via It\^o's formula, from which it is then readily seen that $\pi$ must be $\pi_a^{\gbd}$ for a suitable value of $a$.\\\\

 \noindent {\bf Acknowledgements}
 We would like to thank Prof. B.V. Rao whose insightful and probing questions led to the results reported in Section \ref{sec:extrem}.
Research supported in part by  the RTG award (DMS-2134107) from the NSF.  SB was supported in part by the NSF-CAREER award (DMS-2141621).
AB was supported in part by the NSF (DMS-2152577).

\bibliographystyle{plain}
\bibliography{NSFbib,bibl,scaling,refs,atlas_ref}

\vspace{\baselineskip}

\noindent{\scriptsize {\textsc{\noindent S. Banerjee and A. Budhiraja,\newline
Department of Statistics and Operations Research\newline
University of North Carolina\newline
Chapel Hill, NC 27599, USA\newline
email: sayan@email.unc.edu
\newline
email: budhiraj@email.unc.edu
\vspace{\baselineskip} } }}

\end{document}